\documentstyle[12pt,amscd]{amsart}
\newtheorem{thm}{Theorem}[section]
\newtheorem{cor}{Corollary}[section]
\newtheorem{lem}{Lemma}[section]
\newtheorem{prop}{Proposition}[section]

\theoremstyle{definition}
\newtheorem{defn}{Definition}[section]
\newtheorem{exmp}{Example}[section]
\newtheorem{prob}{Problem}[section]

\theoremstyle{remark}
\newtheorem{claim}{Claim}[section]
\newtheorem{rem}{Remark}[section]
\begin{document}
\title[Seshadri Constants and free pencils]{Seshadri constants,\\
gonality of space curves\\
and restriction of stable bundles}
\author{Roberto Paoletti}
\address{Department of Mathematics \\
UCLA \\
Los Angeles, CA 90024}
\curraddr{Dipartimento di Matematica, Universit\'a di Pavia,
27100 Pavia, Italy}

\maketitle
\bigskip
\bigskip


\section{\bf {Introduction}}

There exist many situations in algebraic geometry
where the extrinsic geometry of a variety is reflected
in clear restrictions in the way that it can map
to projective spaces. For example, it is well-known that
the gonality of a smooth plane curve $C$ of degree $d$ is
$d-1$, and that every minimal pencil has the form
$\cal O_C(H-P)$, where $H$
denotes the hyperplane class and $P\in C$.

In fact, there exist to date various statements of this kind
concerning the existence of morphisms from a divisor to
$\bold P^1$. The first general results in this direction
are due to Sommese (\cite{so:amp}) and Serrano (\cite{se:ext}).
Reider (\cite{re:app}) then showed that at least part of
Serrano's results for surfaces can be obtained by use of
vector bundle methods based on the Bogomolov-Gieseker inequality
for semistable vector bundles on a surface.

In \cite{pa}, a generalization of these methods to higher
dimensional varieties is used to obtain the following statement:

\begin{thm} Let $X$ be a smooth projective $n$-fold,
and let $Y\subset X$ be a reduced irreducible divisor. If
$n\ge 3$ assume that $Y$ is ample, and if $n=2$ assume
that $Y^2>0$ (so that in particular it is at least nef).
Let $\phi :Y@>>>\bold P^1$ be a morphism, and let $F$ denote
the numerical class of a fiber.

\noindent
(i) If
$$F\cdot Y^{n-2}<\sqrt {Y^n}-1,$$
then there exists a morphism $\psi :X@>>>\bold P^1$
extending $\phi$. Furthermore, the restriction
$$H^0(X,\psi ^{*}\cal O_{\bold P^1}(1))@>>>
H^0(Y,\phi ^{*}\cal O_{\bold P^1}(1))$$
is injective. In particular, $\psi$ is linearly normal if
$\phi$ is.

\noindent
(ii) If
$$F\cdot Y^{n-2}=\sqrt {Y^n}-1$$
and $Y^n\neq 4$, then either there exists an
extension $\psi :X@>>>\bold P^1$ of $\phi$, or else
we can find an effective divisor
$D$ on $X$ such that
$(D\cdot Y^{n-1})^2=(D^2\cdot Y^{n-2})Y^n$ and
$D\cdot Y^{n-1}=\sqrt {Y^n}$,
and an inclusion
$$\phi ^{*}\cal O_{\bold P^1}(1)\subset \cal O_Y(D).$$
\end{thm}
\bigskip

However, a much less understood range of situations is the one where
$codim(Y)\ge 2$. In some particular cases there are
rather precise statements. In curve theory, in particular,
one has a clear picture of the gonality of Castelnuovo extremal
curves (\cite{acgh}). In even degree,
for example, if $C\subset \bold P^3$ is a smooth complete intersection
of a smooth quadric and a hypersurface of degree $a\ge 2$, the
gonality is attained by restricting to $C$ the rulings on the quadric.
More generally, unpublished work of Lazarsfeld shows that
if $C\subset \bold P^3$ is a smooth complete intersection of type
$(a,b)$, with $a\ge b$, then $gon(C)\ge a(b-1)$. Lazarsfeld's
argument is also based on Bogomolov's instability theorem.
In a somewhat more general
direction, Ciliberto and Lazarsfeld have studied linear series
of low degree on various classes of space curves (\cite{cl}). Their method
is based on the number of conditions imposed by a linear series
on another.

Naturally enough, one is led to investigate
more general situations. We
shall focus on the gonality
of space curves, and then show how the methods
developped apply to other circumstances as well.
In the codimension $1$ case
we have seen that the self intersection of the divisor
governs the numerical constraint on a free pencil on $Y$.
Loosely speaking, in the
higher codimension case a similar role is played by the
Seshadri costant of the curve.
This is defined as follows.
Consider a smooth curve $C\subset \bold P^3$
and denote by
$$f:X_C@>>>\bold P^3$$
the blow up of $\bold P^3$ along $C$, and by
$$E=f^{-1}C$$
the exceptional divisor.
The {\it Seshadri constant} of $C$ is
\begin{center}
$\epsilon (C)=sup\{\eta \in \bold Q|
f^{*}H-\eta E$ is ample$\}$.
\end{center}
This is a very delicate invariant, and it gathers classical
information such as what secants the curve has and
the minimal degree in which powers of $\cal J_C$ are globally
generated.
For example, if $C\subset \bold P^3$ is a complete intersection
of type $(a,b)$, with $a\ge b$, then $\epsilon (C)=\frac 1a$.
More generally, if $C\subset \bold P^3$ is defined
as the zero locus of a regular section of a rank two vector bundle
$\cal E$,
then we have an estimate
$\epsilon (C)\ge \gamma (\cal E)$, where
$\gamma (\cal E)$ is the Seshadri constant of $\cal E$, defined as
\bigskip
\begin{center}
$\gamma (\cal E)=sup\{\frac nm|S^n\cal E^{*}(m)$ is
globally generated$\}$.
\end{center}
\bigskip
\noindent
It is always true that $\epsilon (C)\ge \frac 1d$.
However, the problem of finding general
optimal estimates $\epsilon (C)$ for an
arbitrary curve seems to be a hard one.
Something can be said,
for example, as soon as $C$ can be expressed
as an irreducible component of a complete intersection of smooth
surfaces.

Interest in Seshadri constants, of course, is not new.
In fact, if $Y$ is a subvariety of any projective variety
$X$, one can define in an obvious way the Seshadri constant
of $Y$ with respect to any polarization $H$ on $X$.
Seshadri constants of points, in particular, have received
increasing attention recently, partly in relation to
the quest for Fujita-type results. A differential geometric
interpretation has been given by DeMailly (\cite{de}).
Seshadri constants of points on a surface have been investigated
by Ein and Lazarsfeld (\cite{el}), who have proved the
surprising fact that
they can be bounded away from zero at all but countably many
points of $S$. However, Seshadri constants of higher dimensional
subvarieties have apparently never been put at use.

What a bound on the gonality of a space curve might look
like is suggested by Lazarsfeld's result. In fact, we may
write $a(b-1)=deg(C)-\frac 1{\epsilon (C)}$, so that for a complete
intersection we have the optimal bound
$$gon(C)\ge d-\frac 1{\epsilon (C)}.$$
Keeping the notation above, let us define
$$H_{\eta}=f^{*}H-\eta E$$
and
$$\delta _{\eta}(C)=\eta \cdot deg(N)-d,$$
where $N$ is the normal bundle of $C$.
For example, for a complete intersection of
type $(a,b)$ with $a\ge b$ we have $\delta _{1/a}(C)=b^2$.
$\delta _{\eta}(C)$ has a simple geometric meaning, that we
explain at the end of Chapter 3.
Our result is
\begin{thm} Let $C\subset \bold P^3$ be a smooth
curve of degree $d$ and Seshadri constant $\epsilon (C)$.
Set $\alpha =min \big \{1,\sqrt d\big (1-\epsilon (C)\sqrt d\big
)\big \}$.
Then
$$gon(C)\ge min \Big \{\frac {\delta _{\epsilon (C)}
(C)}{4\epsilon (C)},
\alpha \Big (d-\frac {\alpha}{\epsilon (C)}\Big )\Big  \}.$$
\end{thm}
This reproduces Lazarsfeld's result if $a\ge b+3$. As another
example, it says that if $a\gg b$ and $C$ is residual to a line
in a complete intersection of type $(a,b)$, then
$gon(C)=ab-(a+b-2)$ (consider the pencil of planes through the
line).
In view of the above, one would expect the above bound
to hold with $\alpha =1$ always, but I have been unable to prove it.

The idea of the proof is as follows. If $A$ is a minimal pencil
on $C$, and if $\pi :E@>>>C$ is the induced projection, one
can define a rank two vector bundle on $X_C$ by the exactness
of the sequence
$$0@>>>\cal F@>>>H^0(C,A)\otimes \cal O_{X_C}@>>>\pi ^{*}A@>>>0.$$
The numerical assuptions then force $\cal F$ to be Bogomolov unstable
w.r.t. $H_{\epsilon (C)}$ (see $\S0$) and therefore a
maximal destabilizing line bundle
$$\cal O_{X_C}(-D)\subset \cal F$$
comes into the picture.
$D$ and $A$ are related by the inequalities coming from the instability
of $\cal F$, and from this one can show that $deg(A)$ is forced
to satisfy the above bound.

By its general nature, this argument can be applied to the study
of linear series on arbitrary smooth subvarieties of $\bold P^r$.
We will not detail this generalization here.
\bigskip

In another direction, similar methods have been used by Bogomolov
(\cite{bo:78} and \cite{bo:svb})
to study the behaviour of a stable bundle on a surface under
restriction to a curve $C$ that is linearly equivalent to
a multiple of the polarization at hand.
For example, it follows from Bogomolov's theorem that
if $S$ is a smooth surface with $Pic(S)\simeq \bold Z$
and $\cal E$ is a stable rank two vector bundle on $S$,
then $\cal E|_C$ is also stable, for every irreducible curve
$C\subset S$ such that $C^2>4c_2(\cal E)^2$.
A more complicated statement holds for arbitrary surfaces.
One can see, in fact, that
this result implies a similar one for surfaces in $\bold P^3$.

In the spirit of the above discussion, one is then led to
consider the problem of the behaviour under restriction to
subvarieties of higher codimension. The inspiring idea,
suggested by the divisor case, should be that when some suitable
invariants, describing some form of "positivity" of the
subvariety, become large with respect to the invariants of
the vector bundle, then stability is preserved under restriction.
Furthermore, if in the divisor case one needs the hypothesis
that $\cal E$ be $\cal O_S(C)$-stable, in the higher codimension
case one should still expect some
measure of the relation  between the geometry of
the subvariety and the stability of the vector bundle to play
a role in the solution to the problem.

In fact, in the case of space curves
the same kind of argument that proves the theorem
about gonality
can be applied to this question. Before explaining the result,
we need the following definition. Recall
that if $X$ is a smooth projective
threefold, $\cal F$ is a vector bundle on $X$ and $L$ and
$H$ are two
nef line bundles on $X$, $\cal F$ is
said to be $(H,L)$-stable if for every
nontrivial subsheaf $\cal G\subset
\cal F$ we have $(fc_1(\cal G)-gc_1(\cal F))\cdot H\cdot L<
0$,
where $f=rank(\cal F)$ and $g=rank(\cal G)$.
Let then $\cal E$ be a rank two vector bundle on
$\bold P^3$, and consider a curve $C\subset \bold P^3$.
Let us define the
{\it stability constant of $\cal E$ with respect to $C$} as
\bigskip
\begin{center}
$\gamma (C,\cal E)=sup\{\eta \in [0,\epsilon (C)]|
f^{*}\cal E$ is $(H,H_{\eta})$-stable$\}.$
\end{center}
\bigskip
\noindent
For example, if $C$ is a complete intersection of type $(a,b)$
and the restriction of $\cal E$ to one of the two surfaces defining $C$
is stable (with respect to the hyperplane bundle) then
$\gamma (C,\cal E)=\epsilon (C)$.

Then we have

\begin{thm} Let
$\cal E$ be a stable rank two vector bundle on $\bold P^3$
with $c_1(\cal E)=0$. Let $C\subset \bold P^3$ be a smooth curve
of degree $d$ and Seshadri constant $\epsilon (C)$,
and let $\gamma =\gamma (C,\cal E)$ be the stability constant
of $\cal E$ w.r.t. $C$. Suppose that $\cal E|_C$ is not
stable. Then
$$c_2(\cal E)\ge min\Big \{\frac {\delta _{\gamma}(C)}
4, \alpha \gamma \Big (d-\frac {\alpha}{\gamma}\Big ) \Big \},$$
where $\alpha =:min\Big \{1,\sqrt d
\Big (\sqrt {\frac 34}-\gamma \sqrt d
\Big )\Big \}$.
\end{thm}

The problem of the behaviour of stable bundles on $\bold P^r$
under restriction to curves has been studied by many researchers.
In particular, a well-known fundamental theorem of Mehta and
Ramanathan
(\cite{mr:res}) shows that $\cal E|_C$ is stable if $C$ is a {\it general}
complete intersection curve of type $(a_1,a_2,\cdots)$, and all the
$a_i\gg 0$. Flenner (\cite{fl}) has then given an explicit
bound on the $a_i$s in term of the invariants of $\cal E$
which makes the conclusion of Mehta and Ramanathan's Theorem true.
On the other hand, here we give numerical conditions that
imply stability for $\cal E|_C$, with no generality assumption
and without restricting $C$ to be a complete intersection.

We have the following applications:
\begin{cor} Let $\cal E$ be a stable rank two vector bundle
on $\bold P^3$ with $c_1(\cal E)=0$ and $c_2(\cal E)=c_2$.
Suppose that $b\ge c_2+2$. If $V\subset \bold P^3$ is a smooth
surface of degree $b$, then $\cal E|_V$ is
$\cal O_V(H)$-stable.
\end{cor}

\begin{cor} Let $\cal E$ be a stable bundle on $\bold P^3$
with $c_1(\cal E)=0$ and $c_2(\cal E)\neq 1$. Suppose that $C=
V_a\cap V_b\subset \bold P^3$ is an
irreducible smooth complete intersection
curve and that $V_a$ is smooth.
Assume furthermore that $a\ge \frac 43b+\frac {10}3$ that
and that $b\ge c_2(\cal E)+2$. Then $\cal E|_C$
is stable.
\end{cor}

\begin{cor} Let $c_2\ge 0$ be an integer and let
$\cal M(0,c_2)$ denote the moduli space of stable rank two vector
bundles on $\bold P^3$. If $a\gg b\gg c_2$ and
$C\subset \bold P^3$ is an irreducible smooth complete intersection of
type $(a,b)$, then
$\cal
M(0,c_2)$ embeds in the moduli space of stable vector
bundles of degree $0$ on $C$.
\end{cor}
\bigskip
\bigskip
This paper covers part
of the content of my Phd thesis at UCLA.
I want to thank Robert Lazarsfeld, my advisor, for introducing me to
Algebraic Geometry and taking continuous interest in my progress.

I am also endebted to a number of people for valuable comments
and discussions; among them, O. Garcia-Prada, D. Gieseker, M. Green,
J. Li and A. Moriwaki.

\section{\bf {Preliminaries}}\label{section:preliminaries}

In this section we state some results that will be
used in the sequel.
The following fact is well-known:
\begin{lem} Let $X$ be a smooth projective variety
and let $Y\subset X$ be a divisor.
Suppose that we have an exact sequence:
$$0@>>>\cal F@>>>\cal E@>>>A@>>>0,$$
where $A$ is a line bundle on $Y$ and $\cal E$ is a rank
two vector bundle on $X$. Let $[Y]\in A^1(X)$ be the divisor
class of $Y$ and let $[A]\in A^2(X)$ be the image of
the divisor class of $A$ under the push forward
$A^1(Y)@>>>A^2(X)$.
Then $\cal F$ is a rank two vector bundle on $X$,
having Chern classes
$c_1(\cal F)=c_1(\cal E)-[Y]$ and $c_2(\cal F)=c_2(\cal E)+[A]
-Y\cdot c_1(\cal E)$.
\label{lem:eltr}
\end{lem}
{\it Proof.} The first statement follows by considering
local trivializations.
As to
the Chern classes of
$\cal F$, we could prove the statement by directly
computing
$$c_t(\cal F)=c_t(\cal E)\cdot c_t(A)^{-1}.$$
However, the following shorter argument proves
that the above equalities hold numerically,
after multiplying both sides by $n-2$ nef divisor classes
(which is what we need).
First of all, the morphism $\cal F@>>>\cal E$ drops rank
along $Y$, and therefore $c_1(\cal F)=c_1(\cal E)-Y$. Let
us consider the second equality.
If $X$ is a surface, the proof is reduced to
a Riemann-Roch computation. If
$dim(X)=3$, let $H$ be any very
ample divisor on $X$, and let $S\in |H|$ be a general smooth surface.
By generality, we may assume that $C=S\cap Y$ is a smooth irreducible
curve.
Then by restriction we obtain an exact sequence on
$S$:
$0@>>>\cal F|_S@>>>\cal E|_S@>>>A|_C@>>>0$.
By applying the statement for the surface case,
we then obtain $(c_2(\cal F)-c_2(\cal E)-[A]
+Y\cdot c_1(\cal E))\cdot H=0$.
But then the expression between brackets has to be
killed by all ample divisors, and so it is numerically
trivial.
The general case is similar.
$\sharp$

\bigskip
\begin{lem} Let $X$ be a smooth projective threefold, and
let $C\subset X$ be a smooth curve in $X$.
Denote by $f:X_C@>>>X$ the blow up of $X$ along $C$, and
let $E$ be the exceptional divisor.
Then $E^3=-deg(N)$, where $N$ is the normal bundle of
$C$ in $X$.
Furthermore, let $A$ be any line bundle on $X$, and by abuse
of language let $A$ also denote its pull-back to $X_C$.
Then $E^2\cdot A=-C\cdot A$.
\label{lem:segre}
\end{lem}
{\it Proof.} Both statements follow from a simple Segre class
computation (see for example \cite{fu}).
$\sharp$
\bigskip

We now recall some known results about instability of
rank two vector
bundles on projective manifolds, which are one of the main tools
in the following analysis. Recall the following notation.
\begin{defn}
If $S$ is a smooth projective surface, $N(S)$ is the vector
space of the numerical equivalence classes of divisors in
$S$; $K^{+}(S)
\subset N(S)$ is the (positive) cone spanned by those
divisors $D$ such that $D^2>0$ and $D\cdot H>0$ for some
polarization on $S$. In general, if $X$ is a smooth
projective $n$-fold and $H$ is a polarization
on it, we shall denote by $K^{+}(X,H)$ the cone of all
numerical classes $D$ in $N(X)$ such that
$D^2\cdot H^{n-2}>0$ and $D\cdot H^{n-1}>0$ (or,
equivalently, $D\cdot R\cdot H^{n-2}>0$ for
any other polarization $R$ on $X$).
\label{defn:poscone}
\end{defn}

\begin{defn} Let $X$ be a smooth projective $n$-fold, and let
$\cal E$ be a rank two vector bundle on $X$, with Chern
classes $c_1(\cal E)$ and $c_2(\cal E)$.
The {\it discriminant} $\Delta (\cal E)\in A^2(X)$ is
$$\Delta (\cal E)=c_1(\cal E)^2-4c_2(\cal E).$$
\label{defn:discr}
\end{defn}

\begin{lem} Let $X$ be a smooth projective
$n$-fold, and let $\cal E$ be a rank two vector bundle
on $X$.
Fix a polarization $H$ on $X$.
Suppose that $\cal L_1, \cal L_2\subset \cal E$
are two line bundles in $\cal E$. Let us
denote by $l_1$ and $l_2$ their $H$-degrees, respectively
(i.e., $l_i=\cal L_i\cdot H^{n-1}$) and let
$e=deg_H(\cal E)=\wedge ^2\cal E\cdot H^{n-1}$ be the $H$-degree
of $\cal E$. Suppose that
$2l_i>e$ for $i=1$ and $i=2$ (in other words,
$\cal L_1$ and $\cal L_2$ make $\cal E$ $H$-unstable).
If $\cal L_2$ is saturated
in $\cal E$, then $\cal L_1\subseteq \cal L_2$.
\label{lem:contains}
\end{lem}
{\it Proof.} Set
$l=min \{l_1,l_2\}$. By assumption, we have
$$2l-e>0.$$
\begin{claim}
If the statement is false, the morphism of vector
bundles
$$\phi :\cal L_1\oplus \cal L_2@>>>\cal E$$
is generically surjective.
\end{claim}
{\it Proof} Set $\cal Q=\cal E/\cal L_2$. Then
$\cal Q$ is a rank one torsion free sheaf. The morphism
$\cal L_1@>>>\cal Q$ is therefore either identically
zero or generically nonzero. If $\cal L_1 \not \subset
\cal L_2$ the morphism $\cal L_1@>>>\cal Q$ is then
generically nonzero. But this implies that
$\phi$ is generically surjective.
$\sharp$
\bigskip

Therefore,
$\wedge ^2\cal E\otimes \cal L_1^{-1}\otimes \cal L_2
^{-1}$ is an effective line
bundle; it follows that
$$0\le e-(l_1+l_2)\le e-2l,$$
a contradiction.
$\sharp$
\bigskip

\begin{cor} Let $X$ and $\cal E$ be as above, and
let $\cal A\subset \cal E$ be a saturated
$H$-destabilizing line bundle.
Then
$\cal A$ is the maximal $H$-destabilizing line bundle.
\label{cor:max}
\end{cor}

\begin{cor}  Let $X$ be a smooth projective $n$-fold,
and fix a very ample
linear series $|V|$ on $X$, with $V\subset
H^0(X,H)$.
Suppose that
$\cal E$ is a rank two vector bundle on $X$ which is
$H$-unstable. Let $C\subset X$ be a general complete intersection
of $n-1$ divisors in $|V|$.
Then the maximal destabilizing line bundle of $\cal E|_C$ is the
restriction to $C$ of the maximal destabilzing line bundle of $\cal
E$.
\label{cor:res}
\end{cor}
{\it Proof.} Let $\cal A$ be the maximal destabilizing line bundle
of $\cal E$. Then the inclusion $
\psi :\cal A@>>>\cal E$ drops rank
in codimension two, because $\cal A$ is saturated in $\cal E$.
Let $Z$ be the locus where $\psi$ drops rank. For a general complete
intersection curve, we have $C\cap Z=\emptyset$.
Hence $\cal A|_C$ is the maximal destabilizing line bundle of
$\cal E|_C$.
$\sharp$
\bigskip

The basic result is the following
\begin{thm} (Bogomolov) Let
$S$ be a smooth projective surface, and let
$\cal E$ be a rank two vector bundle on $S$. Let $c_1(\cal E)$
and $c_2(\cal E)$ be its Chern classes, and suppose that
$$c_1(\cal E)^2-4c_2(\cal E)>0.$$
Then there exists an exact sequence
$$0@>>>A@>>>\cal E@>>>\cal J_Z\otimes B@>>>0,$$
where $A$ and $B$ are line bundles on $S$ and $Z$ is a codimension
two (possibly empty) local complete intersection subscheme,
with the property that
$A-B\in K^{+}(S)$.
\label{thm:bog}
\end{thm}
For a proof, see \cite{bo:st}, \cite{mi:cc}, \cite{re:vbls},
\cite{gi} or \cite{la:svbt}.

\begin{cor} Let $S$ and $\cal E$ be a smooth projective surface
and a rank two vector bundle on it such that
the hypothesis of the theorem are satisfied.
Let $\cal A$ and
$\cal B$ be the line bundles in the above exact sequence.
Then the following inequalities hold:
$$(\cal A-\cal B)\cdot H>0$$
for all polarizations $H$ on $S$, and
$$(\cal A-\cal B)^2\ge c_1(\cal E)^2-4c_2(\cal E).$$
\label{cor:devissage}
\end{cor}
{\it Proof.} The first inequality follows from the condition
$A-B\in K^{+}(S)$.
To obtain the second, just use the above exact sequence to
compute $c_1(\cal E)$ and $c_2(\cal E)$: we obtain
$$c_1(\cal E)^2-4c_2(\cal E)=(A+B)^2-4A\cdot B-4deg[Z]
\le  (A-B)^2.$$
$\sharp$

\bigskip

\begin{cor} Let $S$ and $\cal E$ satisfy the hypothesis of
Bogomolov's theorem, and let $H$ be any polarization on $S$.
Then $\cal E$ is $H$-unstable, and $\cal A$ is the
maximal $H$-destabilizing subsheaf of $\cal E$.
\end{cor}

Recall the fundamental theorem
of Mumford-Mehta-Ramanatan (cfr \cite{mi:cc}):

\begin{thm} Let $X$ be a smooth
projective $n$-fold,and let $H$ be a polarization
on $X$. Consider a vector bundle $\cal E$ on $X$.
If $m\gg 0$ and $V\in |mH|$ is general,
then the maximal $H|_V$-destabilizing subsheaf of
$\cal E|_V$ is the restriction of the maximal $H$-destabilizing
subsheaf of $\cal E$.
\label{thm:mumera}
\end{thm}

This theorem is very powerful, because it detects global instability
from instability on the general complete intersection curve.

\begin{thm} Let $X$ be a smooth projective $n$-fold, and let
$H$ be a polarization on $X$.   Consider a rank two
vector bundle
$\cal E$ on $X$, and suppose that
$$(c_1(\cal E)^2-4c_2(\cal E))\cdot H^{n-2}>0.$$
Then there exists an exact sequence
$$0@>>>\cal A@>>>\cal E@>>>\cal B\otimes \cal J_Z@>>>0,$$
where $\cal A$ and $\cal B$ are invertible sheaves and
$Z$ is a locally complete intersection of codimension two
(possibly empty) such that
$$\cal A-\cal B\in K^{+}(X,H)  $$
and
$$(\cal A-\cal B)^2\cdot H^{n-2}\ge (c_1(\cal E)^2-4c_2(\cal E))
\cdot H^{n-2}.$$
Furthermore, $\cal A$ is the maximal $(H,\cdots ,H,L)$-destabilizing
subsheaf of $\cal E$, for every ample line bundle $L$
on $X$.
\label{thm:main}
\end{thm}

{\it Proof.}
The case $n=2$ is just
the content of Theorem \ref{thm:bog}; for $n\ge 3$, the
statement follows by induction using theorem \ref{thm:mumera}.
$\sharp$
\bigskip

\begin{defn} If $\cal E$ satisfies the hypothesis of the
theorem, we shall say that $\cal E$ is {\it Bogomolov-unstable with
respect to $H$}.
\label{defn:bogunst}
\end{defn}
\bigskip

\begin{lem} Let $f:X@>>>Y$ be a morphism of
projective varieties.
Let $\cal F$
and $A$ be, respectively, a vector bundle and an ample line
bundle on $X$.
For $y\in Y$, let $X_y=f^{-1}y$ and denote by
$\cal J_{X_y}$ the ideal sheaf of $X_y$.
Then there exists $k>0$ such that
$$H^i(X,\cal F\otimes A^n\otimes \cal J_{X_y})=0$$
for all $i>0$, $n\ge k$ and for all $y\in Y$.
\label{lem:fs}
\end{lem}
{\it Proof.} For all $y\in Y$, there is an exact sequence
$$0@>>>\cal F\otimes A^n\otimes \cal J_{X_y}@>>>
\cal F\otimes A^n@>>>\cal F\otimes A^n|_{X_y}@>>>0.$$
Furthermore, there exists a {\it flattening stratification}
of $Y$ w.r.t. $f$, $Y=
\coprod _{l=1}^r Y_l$, with the following property (\cite{mu:cs}).
The $Y_l$ are locally closed subschemes of $Y$, and
if $X_l=:f^{-1}Y_l$, $l=1,\cdots,r$, and $f_l:X_l@>>>Y_l$
is the restriction of $f$, then $f_l$ is a flat morphism.
Let us then start by finding $k_1$ such that
for all $n\ge k_1$ we have
$$H^i(X,\cal F\otimes A^n )=0$$
and
$$H^i(X,\cal F\otimes A^n\otimes \cal J_{X_l})=0$$
for all $i>0$ and for all $l=1,\cdots,r$.
Then
it is easy to see that the statement is equivalent to saying that there is
$k\ge k_1$ such that for all $n\ge k$ the restriction maps
$$H^0(X,\cal F\otimes A^n)@>\phi _y>>H^0(X_y,\cal F\otimes
A^n|_{X_y})$$
are all surjective, and that
$$H^i(X_y,\cal F\otimes A^n|_{X_y})=0,$$
for all $y\in Y$ and for all $i>0$.
If $y\in Y_l$ and $\cal J^{X_l}_{X_y}$ denotes the
ideal sheaf of $X_y$ in $X_l$, then we have an exact sequence
$$0@>>>\cal J_{X_l}@>>>\cal J_{X_y}@>>>\cal J^{X_l}_{X_y}@>>>0.$$
\begin{claim} The lemma is true if there exists $k$ such that
for all $n\ge k$, for $l=1,\cdots,r$ and for all $y\in Y_l$ we
have that $H^i(X_l,\cal F\otimes A^n|_{X_l}\otimes
\cal J^{X_l}_{X_y})=0$ for $i>0$.
\end{claim}
{\it Proof.} It follows from the exact sequences
$$
\CD
H^i(X,\cal F\otimes A^n\otimes \cal J_{X_l})@>>>H^i(X,\cal F\otimes A^n
\otimes \cal J_{X_y})
@>>>H^i(X,\cal F\otimes A^n\otimes \cal J^{X_l}_{X_y})\\
@|                   @.              @|   \\
0   @.    @.        0
\endCD
$$
for $i>0$.
$\sharp$
\bigskip

This means that we can reduce to the case where $f$ is flat.
For $y_0\in Y$, we can find $k_0$ such that
for $n\ge k_0$ and for $i>0$ we have
$$H^i(X,\cal F\otimes A^n\otimes \cal J_{X_{y_0}})=0.$$
Therefore, the morphism
$$\lim_{y_0\in U}H^0(f^{-1}U,\cal F\otimes A^n)@>\beta _{y_0}>>
H^0(X_y,\cal F\otimes A^n)
$$
is onto, and then so is
$$\psi _{y_0}=:\beta_{y_0} \otimes k(y_0):
f_{*}(\cal F\otimes A^n)(y)@>>>H^0(X_{y_0},\cal F\otimes A^n
|_{X_{y_0}}).$$
By Grauert's theorem (\cite{ha:ag}) we then have that $\psi _{y_0}$
is an isomorphism, and that the same holds for
$\psi _y$, for $y$ in a suitable open neighbourhood $U_0$
of $y_0$.
Therefore the restriction morphism
$$H^0(X,\cal F\otimes A^n)@>>>H^0(X_y,\cal F\otimes A^n|_{X_y})$$
come from a morphism of sheaves, and hence
they are onto for all $y\in V_0$,
for a suitable open set $V_0\subset U_0$.
We can then invoke the quasi-compactness of $Y$ to
conclude that there exists $k$ such that
$H^1(X,\cal F\otimes A^n\otimes \cal J_{X_y})=0$ for all
$y\in Y$.
As to $i\ge 2$, we have isomorphisms
$$H^i(X_y,\cal F\otimes A^n|_{X_y})
\simeq H^{i+1}(X,\cal F\otimes A^n\otimes \cal J_{X_y})=0$$
for all $i>0$, and so we need to show that
$H^i(X_y,\cal F\otimes A^n)=0$ for $n\gg 0$,
$i>0$ and for all $y\in Y$.
But for $n\gg 0$ we have
$$R^if_{*}(\cal F\otimes A^n)=0$$
if $i>0$ and then this implies $h^i(X_y,\cal F\otimes A^n|_{X_y})=0$
for all $y$ (\cite{mu:cs}).
$\sharp$
\bigskip

We record here a trivial numerical lemma that will be handy in
the sequel:
\begin{lem} If $s\ge \alpha$, $a \ge 2s$ and
$b \ge a s-s^2$, then $b \ge a \alpha-\alpha ^2$.
\label{lem:trivial}
\end{lem}
{\it Proof.} $a s-s^2$ is increasing
in $s$ if $a \ge 2s$.
The statement follows.
$\sharp$
\bigskip

\section{\bf {Seshadri Constants of Curves}}\label{section:sc}

Let $C\subset \bold P^3$ be a smooth curve
and
let $H$ denote the hyperplane bundle on $\bold P^3$.
We shall let $f:X_C@>>>\bold P^3$ be the blow up of $\bold P^3$ along $C$,
and $E=f^{-1}C$ be the exceptional divisor.
\begin{defn} The Seshadri constant of $C$ is
\begin{center}
$\epsilon (C)=:
sup \{\eta \in \bold Q|f^{*}H-\eta E$ is ample$\}$.
\end{center}
\label{defn:sc}
\end{defn}
\noindent
In other terms, $\epsilon (C)$ is the supremum of the ratios
$\frac nm$, where $n$ and $m$ are such that
$mH-nE$ is ample (or, equivalently, very ample).
In the sequel we shall use the short hand
$$H_{\eta}=:H-\eta E$$
for $\eta \in \bold Q$; furthermore, we shall generally write
$H$ for $f^{*}H$ (as we just did).

\begin{lem} $H_{\eta}$ is ample if and only if
$0<\eta <\epsilon (C)$. It is nef if and only if
$\eta \in [0,\epsilon (C)]$.
\label{lem:nef}
\end{lem}
{\it Proof.} Since the ample cone of a projective
variety is convex, the line $H-tE\subset N^1(X)$ intersects
$K^{+}(X)$ in a segment $(H-t_1E,H-t_2E)$.
Let $F$ denote the numerical class of a fiber
of $\pi :E@>>>C$. Then $H_{\eta}\cdot F=\eta$, and therefore if
$H_{\eta}$ is ample we must have $\eta >0$.
Hence $t_1\ge 0$. On the other hand, it is well known
that $H-tE$ is ample for $t>0$ sufficiently small,
and therefore $t_1=0$.
By definition,
$t_2=\epsilon _2(C)$.
The remaining part of the statement is clear.   $\sharp$
\bigskip

\begin{cor} We have
\begin{center}

$\epsilon (C)=sup\{\eta |H_{\eta}\cdot D\ge 0$ for all curves
$D\subset X_C\}.$
\label{cor:nef}
\end{center}
\end{cor}

\begin{lem} Let $C\subset \bold P^3$ be a smooth curve,
and let $\cal J_C$ be its ideal sheaf.
Let $m$ and $n$ be nonnegative integers.
Then $\cal O_{X_C}(mH-nE)$
is globally generated if $\cal J_C^n(m)$ is.
\label{lem:gg}
\end{lem}
{\it Proof.} Let us suppose that $\cal J_C^n(m)$ is globally generated,
and let
$$F_1,\cdots,F_k\in H^0(\bold P^3,\cal J_C^n(m))$$
be a basis.
Let $P\in C$ and
let $U$ be some open
neighbourhood of $P$.
By assumption, $F_1,\cdots,F_k$ generate $\cal J_C$ in
$U$. By abuse of language, let us write
$F_i$ for the pull-backs $f^{*}F_i$. Then if $e$ is a local
equation for $E$ in a Zariski open set
$V\subset f^{-1}U$, then the ideal generated by the $F_i$s
is $<\{F_i\}>=(e^n)$.
Hence we can write
$$\sum _{i=1}^kP_iF_i=e^n$$
for some $P_i$s regular on $V$. However, by construction we
can write $F_i=\tilde F_ie^n$, and therefore
we have
$$\sum _{i=1}^k\tilde F_iP_i=1$$
in $V$. Hence the $\tilde F_i$ are base point free, and
they can be extended to global sections of
$\cal O_{X_C}(mH-nE)$, which is therefore globally spanned.
$\sharp$

\bigskip

\begin{cor} Let $C\subset \bold P^3$
be a smooth curve.
Then
\begin{center}
$\epsilon (C)\ge sup \{\frac nm|\cal J_Y^n(m)$ is globally
generated$\}$.
\end{center}
\label{cor:gg}
\end{cor}

Let us look at some examples.
\begin{exmp} If $L\subset \bold P^3$ is
a line, then $\cal J_L(1)$ is globally generated.
Therefore $\epsilon (L)\ge 1$.
On the other hand, let $\Lambda \subset \bold P^3$ be a hyperplane
containing $L$ and let $D\subset \Lambda$ be any irreducible curve
distinct from $L$. Then $H_1\cdot \tilde D=deg(D)-L\cdot _{\Lambda}D=0$,
where $\tilde D\subset Bl_L(\bold P^3)$ is the proper transform of $L$.
Hence $\epsilon (L)=1$. As we shall see shortly, this generalizes
to the statement that if $C\subset \bold P^3$ is a smooth complete
intersection of type $(a,b)$ and $a\ge b$, then
$\epsilon (C)=\frac 1a$.
\end{exmp}
\begin{exmp} If $C\subset \bold P^3$ has an $l$-secant line, then
$\epsilon (C)\le \frac 1l$. To see this, let $L$ be the $l$-secant;
denoting by $\tilde L\subset X_C$ the proper transform of $L$
in $Bl_C(\bold P^3)$ we have $H\cdot \tilde L=1$ and $E\cdot \tilde L=l$.
Hence $0\le H_{\epsilon}\cdot \tilde L$ implies $\epsilon \le \frac 1l$.
\end{exmp}

\begin{lem} Let $C\subset \bold P^3$ be a smooth
curve of degree $d$. Then
$$\frac 1{\sqrt d}\ge \epsilon (C)\ge \frac 1d$$
\label{lem:deg}
\end{lem}
{\it Proof.} It is well-known that a
smooth subvariety of degree $d$ of
projective space is cut out by hypersurfaces of degree $d$.
Hence $\cal J_C(d)$ is globally generated, and this proves the
second inequality.
As to the first, we must have $0\le H\cdot H_{\epsilon}^2=1-\epsilon ^2d$,
by a simple Segre class computation.
$\sharp$
\bigskip

The right inequality is sharp if the curve is degenerate; the left
one is sharp for a complete intersection curve of type $(a,a)$.
If the curve is nondegenerate, however, one can say something
more.
\begin{defn} Let $C\subset \bold P^3$ be a smooth
curve, and let $\cal J_C$ be its ideal sheaf.
$C$ is said to be {\it $l$-regular} if
$$H^i(\bold P^3,\cal J_C(l-i))=0$$
for all $i>0$. The {\it regularity} of $C$, denoted by $m(C)$,
is the smallest $l$ such that $C$ is $l$-regular
(\cite{ca}, \cite{mu}, \cite{gr}).
\label{defn:reg}
\end{defn}
\begin{rem} By a celebrated theorem
of Castelnuovo, we
have $m(C)\le d-1$ (\cite{ca}, \cite{glp}).
\label{rem:cast}
\end{rem}

\begin{prop} Let $C
\subset \bold P^3$ be a smooth space curve, and let $m=m(C)$
be its regularity. Then
$$\dfrac 2{m-1}\ge \epsilon (C)\ge \dfrac 1m.$$
\label{prop:reg}
\end{prop}
{\it Proof.} By a classical theorem of
Castelnuovo-Mumford, the homogeneuos ideal of $C$ is saturated in
degree $m(C)$ and therefore $\epsilon (C)\ge \frac 1{m(C)}$.
By definition, to prove the first inequality it is
enough to show that
$H^i(\bold P^3,\cal J_C(k))=0$ for $k\ge \lceil \frac 2{\epsilon (C)} \rceil
-3$
because this implies $m(C)\le \frac 2{\epsilon (C)}+1$
and then the statement.
To prove the above vanishing, observe that
$$\Big \{2/ (\lceil   2/\epsilon (C)\rceil +1)
\Big \} <\epsilon (C)$$
and therefore
$$\Big (\Big \lceil \frac 2{\epsilon (C)}\Big \rceil +1\Big )H-2E$$
is an ample integral divisor in $X_C$.
Since $\omega _{X_C}=\cal O_{X_C}(-4H+E)$,
the Kodaira vanishing theorem gives:
$$H^i(X_C,\cal O_{X_C}((\lceil  2/\epsilon (C) \rceil
-3)H-E))=0$$
for $i>0$, as desired.
$\sharp$
\bigskip

\begin{rem} Using vanishing theorems on the blow up to
obtain bounds on the regularity is a well-known technique:
see \cite{bel} for various results in this direction.
\end{rem}

\begin{rem} It is not possible, in
the above vanishing, to replace the
condition on $k$ by $k\ge \lceil \frac 1{\epsilon}\rceil$. To see this,
suppose that $C$ is a complete intersection of type $(a,b)$ so that we
have a Koszul resolution
$$0@>>>\cal O_{\bold P^3}(-b)@>>> \cal O_{\bold P^3}\oplus
\cal O_{\bold P^3}(a-b)@>>>\cal J_C(a)@>>>0.$$
It follows that $H^2(\bold P^3,\cal J_C(a))\simeq H^3(\bold P^3,
\cal O_{\bold P^3}(-b))\neq 0$, for $b\ge 4$.
\end{rem}

\begin{cor}  Let $C\subset \bold P^3$ be a
nondegenerate smooth curve. Then
$\epsilon (C)\ge \frac 1{d-1}$.
\end{cor}

\noindent
Equality is attained in the previous corollary in the case of a twisted
cubic.
\bigskip

It is convenient to introduce the following definition.
\begin{defn} Let $C\subset \bold P^3$ be a smooth curve.
For an irreducible curve $D\subset \bold P^3$ different from $C$
let $\tilde D$ be its proper transform in the blow up of $\bold P^3$
along $C$.
Define
\begin{center}
$\epsilon _1(C)=:sup \{\eta \in \bold Q|
(H-\eta E)|_E$ is ample$\}$
\end{center}
and
\begin{center}
$\epsilon _2(C)=:sup\{\eta \in \bold Q|
H_{\eta}\cdot \tilde D\ge 0 \forall$ irreducible curves $D\neq  C$\}.
\end{center}
\label{defn:12}
\end{defn}
\begin{rem} $\epsilon (C)=min\{\epsilon _1(C),
\epsilon _2(C)\}$.
\label{rem:12}
\end{rem}

We are interested in estimating the Seshadri constant of a space curve $C$.
It is convenient to examine $\epsilon _1(C)$ and
$\epsilon _2(C)$ separately. We shall see that $\epsilon _1(C)$ is
determined by the structure of the normal bundle, while $\epsilon _2(C)$
depends on the "linkage" of $C$, and is generally much harder
to estimate.
We start with an analysis of $\epsilon _1(C)$.

\begin{defn} Let $C$ be a smooth projective curve and let
$\cal E$ be a rank two vector bundle on it. For all finite morphisms
$f:\tilde C@>>>C$ and all exact sequences of locally free shaves on $\tilde
C$ of the form $0@>>>L@>>>f^{*}\cal E@>>>M@>>>0$, consider the
ratios $\frac {deg(L)}{deg(f)}$. Let $\Sigma _{\cal E}$ denote the
set of all the numbers obtained in this way. Define
$$s(\cal E)=:sup\Sigma _{\cal E}.$$
\label{defn:wahl}
\end{defn}

\begin{rem}
As in \cite{w}, $s(\cal E)$ can be interpreted as a measure
of the instability of $\cal E$. More precisely, we have
$$s(\cal E)=\frac 12deg(\cal E)$$ if
$\cal E$ is semistable
and
$$s(\cal E)=deg(L)$$
if $\cal  E$ is unstable, and
$L\subset \cal E$ is the maximal destabilzing line subundle
of $\cal E$.
In other words, $s(\cal E)-\frac 12deg(\cal E)\ge 0$ always,
and equality holds if and only if $\cal E$ is semistable.
\label{rem:wahl}
\end{rem}

We then have
\begin{prop} Let $C\subset \bold P^3$ be a smooth
curve.
Denote by $N$ the normal bundle of $C$ in $\bold P^3$, and
let $\epsilon _1(C)$ be as above. Then
$$\epsilon _1(C)=\frac {deg(C)}{s(N)}.$$
\label{prop:e1}
\end{prop}
{\it Proof.} Let $X_C=:Bl_C(X)@>f>>X$ be the blow up of $X$ along
$C$ and let $E$ be the exceptional divisor;
recall that $E$ can be identified with the relative projective
space of lines in the vector bundle $N$.
Set $\pi =f|_E$ and
denote by $F$ a fiber of $\pi$.
Let $D\subset E$ be any reduced irreducible curve. If $D$ is a fiber
of $\pi$, then $\eta >0$ ensures that $H_{\eta}\cdot D>0$. Hence we
may assume that $D@>>>C$ is a finite map, whose degree is given
by $a=D\cdot F$.
Let $
\psi :
\tilde D@>>>D\subset X_C$ be the normalization of $D$, and let
$p:\tilde D@>>>C$ be the induced morphism.
Then $\psi$ is equivalent to the assignment of a
sub-line bundle $L\subset p^{*}N$, given by
$L=\psi ^{*}\cal O_{\bold PN}(-1)$.
Since $\cal O_{\bold PN}(-1)\simeq \cal O_E(E)$,
we have $deg(L)=D\cdot E$. Hence
$H_{\eta}\cdot D=aH\cdot C-\eta \cdot deg(L)$; the condition
$\eta \le \epsilon _1(C)$ translates therefore
in the condition $\eta \le inf\{\dfrac {H\cdot C}{deg(L)/a}\}$.
In other words, then, it is equivalent to
$\eta \le \dfrac {H\cdot C}{s(N)}$.
$\sharp$

\bigskip

\begin{exmp} Let $C\subset \bold P^3$ be a smooth
complete intersection curve of type $(a,b)$, with $a\ge b$.
Then we have a Koszul resolution of the ideal sheaf of $C$,
from which it is easy to conclude that
$\epsilon (C)\ge \frac 1a$. On the other hand, $s(N)=a^2b$
and therefore
by Proposition \ref{prop:e1}
$\epsilon _1(C)=\frac 1a$. Hence we have
$\epsilon (C)=\frac 1a$.
\end{exmp}
\begin{exmp} Let $C\subset \bold P^3$ be given as the zero
locus of a regular section of a rank two vector bundle
$\cal E$ on $\bold P^3$. It is well known that this is always
the case provided that the determinant of the normal bundle $N$
extends. The Koszul resolution then is
$$0@>>>det(\cal E)^{-1}@>>> \cal E^{*}@>>>\cal J_C@>>>0.$$
By Corollary \ref{cor:gg}
and Proposition \ref{prop:e1}, we then conclude that
$$\frac {H\cdot C}{s(\cal E|_C)}\ge \epsilon (C)\ge \epsilon (\cal E)$$
where
\begin{center}
$\epsilon (\cal E)=sup \{\frac nm|S^n\cal E^{*}(m)$ is spanned$\}$.
\end{center}
\label{exmp:zl}
\end{exmp}
We shall call $\epsilon (\cal E)$ the Seshadri constant of the
vector bundle $\cal E$. It has the following geometric interpretation.
Let $\bold P\cal E$ be the relative projective
space of lines in $\cal E$. $Pic (\bold P\cal E)$ is generated
by two line bundles $H$ and $\cal O(1)$, where
$H$ is the pull-back of the hyperplane bundle on $\bold P^3$.
Let $R$ be some divisor
associated to the line bundle $\cal O(1)$. It is well known that
the rational divisor $H+\eta R$ is ample, for
sufficiently small $\eta \in \bold Q^{+}$ (\cite{ha:ag}).
\begin{prop} $\epsilon (\cal E)=sup\{\eta \in \bold Q|
H+\eta R\in Div _{\bold Q}(\bold P\cal E)$ is ample$\}$.
\label{prop:scvb}
\end{prop}
{\it Proof.} Provisionally denote by $\gamma (\cal E)$
the right hand side of the statement. Also, for brevity
let us set $X=\bold P\cal E$ and let $X_z$ stand for the
fiber over a point $z\in \bold P^3$.
Let us first prove that $\epsilon (\cal E)\le \gamma (\cal E)$.
Suppose then that $\eta =\frac nm<\epsilon (C)$, where $n$ and
$m$ are such that $S^n\cal E^{*}(m)$ is globally
generated. Since
$$S^n\cal E^{*}(m)=f_{*}\cal O_X(mH+nR),$$
we have the identifications
$$H^0(X,\cal O_X(mH+nR))\simeq H^0(\bold P^3, S^n\cal E^{*}(m))$$
and
$$H^0(X_z,\cal O_{X_z}(mH+nR))\simeq
S^n\cal E^{*}(m)(z).$$
With this in mind, we then have a surjection
$$H^0(X,\cal O_X(mH+nR))@>>>H^0(X_z,\cal O_{X_z}(mH+nR))$$
for all $z\in \bold P^3$, and since $\cal O_X(mH+nR)$ is generated along
the fibers, it is also globally generated.

Let us now prove that $\gamma (\cal E)\le \epsilon (\cal E)$.
Let $\eta =\frac nm <\gamma (\cal E)$, where $n$ and $m$ have been
chosen so that $mH+nR$ is ample. After perhaps multiplying $m$ and $n$
by some large positive integer we may suppose that $mH+nR$ is very ample and
that
$$H^i(X,\cal J_{X_z}(mH+nR))=0$$
for all $i>0$ and all $z\in \bold P^3$ (see Lemma \ref{lem:fs}).
But then we have surjective restriction maps
$$H^0(X,\cal O_X(mH+nR))@>>>H^0(X_z,\cal O_{X_z}(mH+nR))$$
for all $z\in \bold P^3$, and the lemma then follows from the above
identifications.
$\sharp$

\begin{rem}
The inequality $\epsilon (C)\ge \epsilon (\cal E)$
from Example \ref{exmp:zl}
can then be
explained as follows. For each $n\ge 0$ we have surjective
morphisms
$S^n\cal E^{*}@>>>\cal J_C^n$,
and therefore we have a surjection of sheaves of graded algebras
$$\bigoplus _{n\ge 0} S^n\cal E^{*}@>>>
\bigoplus _{n\ge 0}\cal J_C^n,$$
which yields a closed embedding
$$i:X_C\hookrightarrow \bold P\cal E.$$
On the other hand, $i^{*}\cal O_{\bold P\cal E}(R)=\cal O_{X_C}(-E)$
and the above ineqality is just saying that if $H+\eta R$ is
ample, it restricts to an ample divisor on $X_C$.
\end{rem}
\bigskip

We now consider ways to estimate $\epsilon _2(C)$.
$\epsilon _2(C)$ gathers more global information than
$\epsilon _1(C)$, because it relates to how $C$ is
"linked" to the curves in $\bold P^3$.
Recall that our definition was:
\begin{center}
$\epsilon _2(C)=:sup\{\eta \in \bold Q|
H_{\eta}\cdot \tilde D\ge 0$ $\forall$ irreducible curves $D\subset
\bold P^3$, $D\neq C \}$.
\end{center}

As usual, $\tilde D$ denotes the proper transform of
$D$ in the blow up of $C$.
There does not seem
to be much that one can say about $\epsilon _2(C)$
in general; with some extra assumptions, however, we can obtain
an estimate.

Let us make the following definiton:
\begin{defn} Let $D\subset \bold P^3$ be a reduced irreducible curve,
and let $t:D_n@>>>D\subset \bold P^3$ be its normalization.
If the derivative $dt:T_{D_n}@>>>t^{*}T_{\bold P^3}$ never drops rank, we
shall say that $D$ has only ordinary singularities.
\label{defn:os}
\end{defn}

\begin{prop} Let $C\subset \bold P^3$ be a smooth
curve. Suppose that $C$ is contained in the intersection
of two distinct reduced and irreducible hypersurfaces
$V_a$ and $V_b$ of degree $a$ and $b$, respectively.
Suppose
that all the residual curves to $C$ in the complete intersection
$V_a\cap V_b$ are reduced and
that at least one of the two hypersurfaces is smooth.
Then $$ \epsilon _2(C)\ge \frac 1{a+b-2}.$$
If all the residual curves
have ordinary singularities, then
equality holds if and only if the residual curve is a
union of disjoint lines.
\label{prop:e2}
\end{prop}

\begin{exmp} It is well-known that a curve which is
linked to a line
$L$ in a complete intersection of type $(a,b)$
is cut out by the hypersurfaces $V_a$ and $V_b$ and by a third
equation of degree $a+b-2$. Therefore its ideal sheaf is
generated in degree $a+b-2$, so that
$\epsilon (C)\ge \frac 1{a+b-2}$. On the other hand, it is
easy to
check that $\tilde L\cdot E=a+b-2$. Therefore in this
case we find directly that $\epsilon (C)=\frac 1{a+b-2}$.
More generally, the same argument works whenever $C$ is linked
to a union of (reduced) disjoint lines.
\label{exmp:sharp}
\end{exmp}
\begin{exmp} The assumption that the residual curves be
all reduced is necessary. To see this, let $L\subset \bold P^3$
be a line, and let $V$ be a smooth surface of degree $v$ through
$L$.
We have $L\cdot _VL=2-v$. Let $H$ be the hyperplane bundle
restricted to $V$. Then for $s\gg 0$ the linear series
$|sH-2L|$ is very ample. Choose a smooth curve
$C\in |sH-2L|$. Then $C$ is linked to a double line supported
on $L$ in the complete intersection $V\cap W$, where
$W$ is a suitable hypersurface of degree $s$ in $\bold P^3$.
We have
$$\tilde L\cdot E_C=(sH-2L)\cdot _VL=s-2L^2=s+2v-4,$$
and so $\epsilon _2(C)\le \dfrac 1{s+2v-4}$.
\end{exmp}
{\it Proof.}
We need to show that for $\eta \le\frac 1{a+b-2}$ we
have $\tilde D\cdot H_{\eta}\ge 0$, whenever $D\subset \bold P^3$
is some irreducible curve distinct from $C$. Clearly we
may assume that $D$ is reduced.

Let us start with the following simple observation. Let
$C$ and $D$ be reduced curves in $\bold P^3$, and let
$D_n@>t>>D\subset \bold P^3$ be the normalization of
$D$. If $X_C@>f>>\bold P^3$ is the blow up of $C$ and
$E_C$ is the exceptional divisor, clearly $t$ factors
through $f$, i.e. there exists $u:D_n@>>>X_C$ such that
$t=f\circ u$. On the other hand, $t^{-1}C=u^{-1}f^{-1}
C=u^{-1}E_C$ and therefore
\begin{equation}
\tilde D\cdot E_C=D_n\cdot _uE_C=deg\{t^{-1}C\}.
\label{eq:norm}
\end{equation}

Given the geometric situation,
we start testing the
desired positivity condition on the curves that are not
contained in $V_a\cap V_b$.

\begin{lem} Let $C\subset \bold P^3$, $V_a$ and
$V_b$ be as in the statement of the Proposition.
Suppose that $a\ge b$, and let $\eta \le \frac 1a$.
Then for every irreducible curve $D\not\subset V_a\cap V_b$
we have $\tilde D\cdot H_{\eta}\ge 0$.
\label{lem:nonres}
\end{lem}
{\it Proof of the Lemma.} Let $D$
be reduced and have degree $s$, and set
$G=:V_a\cap V_b$.
$G$ is a complete intersection curve, and then we know
from the Koszul resolution of its ideal sheaf that its
Seshadri constant satisfies
$\epsilon (G)\ge \frac 1a$.
Let $X_G@>>>\bold P^3$ be the blow up of $\bold P^3$ along
$G$, and let $E_G$ be the exceptional divisor.
For $\eta \in \bold Q$, let $H^{\prime}_{\eta} =
g^{*}H-\eta E_G$.
By what we have just said, $H^{\prime}_{\frac 1a}$ is
a nef divisor on $X_G$.
Therefore, if we let
$D^{\prime}\subset X_G$ denote the proper transform
of $D$ in $X_G$, we have
$D^{\prime}\cdot H^{\prime}_{\eta}\ge 0$, and this
can be rewritten as $D^{\prime}\cdot E_G\le as$.
Now let as above $t:D_n@>>>D\subset \bold P^3$
be the normalization of $D$, and let
$\tilde D\subset X_C$ denote the proper transform of
$D$ in the blow up of $C$. Then by equation
(\ref{eq:norm})
$$\tilde D\cdot E_C=deg\{t^{-1}C\}\le deg\{t^{-1}G\}
=D^{\prime}\cdot E_G,$$
since $G\supset C$ as schemes. Therefore,
\begin{equation}
H_{\frac 1a}\cdot \tilde D\ge
H^{\prime}_{\frac 1a}\cdot D^{\prime}\ge 0,
\end{equation}
and the statement follows.
$\sharp$

\bigskip
We now need to consider the
condition $H_{\eta}\cdot \tilde D_i\ge 0$,
where the $D_i$s are the irreducible
components of the residual curve
to $C$ in the complete intersection
$V_a\cap V_b$.
Let us drop the index $i$, and let
$D$ be one of the the $D_i$s.
We have to show that $H_{\eta}\cdot \tilde D\ge 0$ for
$\eta \le \frac 1{a+b-2}$.
We shall be using case (b) of the following lemma, but it may
be worthwhile to state it in more generality:
\begin{lem} Let $C$ and $D$ be reduced irreducible space curves,
and suppose that either one of the following
conditions holds:

(a) $C$ and $D$ are both smooth, or

(b) $C$ and $D$ lie in a smooth hypersurface
$S\subset \bold P^3$.

Then $\tilde D\cdot E_C=\tilde C\cdot E_D$, where $\tilde D$
(resp., $\tilde C$) is the proper transform of $D$ in the
blow up of $C$ (resp., the proper transform of $C$ in the blow up of
$D$).
\label{lem:blowups}
\end{lem}
{\it Proof.} Let us first suppose that
$(b)$ holds.
Let $t: C_n@>>>C\subset S$ be the
normalization of $C$. By (\ref{eq:norm}), we know that
$\tilde C\cdot E_D=deg\{t^{-1}D\}=deg\{t^{*}\cal O_S(D)\}
=C\cdot _SD$
and similarly for
$\tilde C\cdot E_D$.

If (a) holds,
the situation is almost the same, because at each intersection
point $P$ of $C$ and $D$ we can still
locally view $C$ and $D$ as lying
in some smooth
open surface in an neighbourhood of $P$, and
the problem is local in $P$.
Explicitly, the argument is the following.
Suppose that $C\cap D$ is supported on $P_1,\cdots,P_k$.
We "measure" the intersection of $C$ and $D$
in the following way (cfr \cite{sev}):
let $\pi:\bold P^3--\to \bold P^2$ be a general
projection, and set
\begin{equation}
C{*}D=:\sum _{i=1}^ki(\pi (P_i),\pi (C),\pi (D)),
\label{eq:int}
\end{equation}
where $i$ denotes the ordinary intersection multiplicity.
Using the projection
formula, one can easily check the following:
\begin{claim} Let $P\in \bold P^3$
be chosen generally, and let $C_P$ be the cone on $C$ with vertex
$P$. Then
$$C{*}D=\sum _{i=1}^k i(P_i,D,C_P).$$
\label{claim:cones}
\end{claim}
Observe that these intersection multiplicities
are generally constant by the principle of
continuity.
Given that $C*D$ is symmetric, Lemma \ref{lem:blowups} will follow once
we establish that
$C*D=\tilde D\cdot E_C$.

Since $C$ is smooth, it is defined scheme-theoretically
by the cones through it (\cite{mu}). Hence for the proof of Lemma 3.5
we are reduced to the following:

\begin{lem} Let $C$ and $D$ be distinct reduced irreducible
curves in $\bold P^3$. Suppose that $C\cap D$ is supported at points
$P_1,\cdots,P_k$.
Let $\cal C\subset H^0(\bold P^3,
\cal J_C(m))$ be an irreducible family of hypersurfaces.
Suppose that the linear series $V=|\cal C|$ spanned by
$\cal C$ globally generates $\cal J_C(m)$ (in other words,
$C$ is cut out scheme-theoretically by the elements of
$\cal C$).
Then for a general $F\in \cal C$ we have
$$\tilde D\cdot E_C=\sum _{i=1}^k i(P_i,D,F).$$
\label{lem:blowups1}
\end{lem}
{\it Proof.} The assumption implies in particular that
$\cal C\not\subset H^0(\bold P^3,\cal J_C^2(m))$, i.e. that
the general $F\in \cal C$ is generically smooth along
$C$.
For such a general $F$, then, if $\tilde F$ denotes the proper
transform in $X_C$ we have
$$\tilde F\in |f^{*}F-E|.$$
Furthermore, the family of all such $\tilde F$ has to be base
point free, so there is
$F\in \cal C$ which is generically smooth along
$C$ and such that
$\tilde F$ does not meet any of the intersection points of
$\tilde D$ and $E_C$.
Let us denote by a subscript $(\cdot ,\cdot )_{P}$
the contribution to a given intersection product on $X_C$ coming from
the points lying over $P\in \bold P^3$.
Then by construction and the projection formula we have
$$(\tilde D\cdot E_C)_{P_i}=
(\tilde D\cdot f^{*}F)_{P_i}=i(P_i,D,F)$$
and this proves the lemma.
$\sharp$
\bigskip

Let then $X_D@>>>\bold P^3$ be the blow up of $\bold P^3$
along $D$, and let $G$ be the complete intersection
$V_a\cap V_b$.
Then $C$ is a component of the effective cycle
$G-D$, and furthermore $G-D$ does not have any component
supported on $D$.
Hence we may consider the proper transform
$\tilde {G-D}\subset X_D$, which is an effective cycle
in $X_D$ containing $\tilde C$ as a component.
Suppose, say, that $V_a$ is smooth. Then
we are in case (b) of lemma 3.5, and
therefore we have
\begin{equation}
\tilde D\cdot E_C=\tilde C\cdot E_D\le (\tilde {G-D})\cdot E_D.
\label{eq:fulton}
\end{equation}
In the hypothesis of the proposition, at a generic point of
$D$ $V_a$ and $V_b$ are both smooth and meet transversally
(for otherwise $D$ would not be reduced).
Therefore $\tilde V_a\equiv f^{*}V_a-E$ and $\tilde V_b\equiv f^{*}V_b
-E$, and no component of $\tilde V_a\cap \tilde V_b$ maps
dominantly to $D$.
Furthermore if, say, $V_a$ is smooth then
$\tilde V_a\simeq V_a$ does not contain any fiber of $\pi$.
Therefore $\tilde {(G-D)}=\tilde V_a\cap \tilde V_b$, and so
$$\tilde {(G-D)}\cdot E_C= (f^{*}V_a-E)\cdot (f^{*}V_b-E)
\cdot E_C.$$
Let $N$ denote the normal bundle to the complete intersection $G$.
{}From intersection theory, the latter term is known to be
\begin{equation}
\{c(N)\cap s(D,\bold P^3)\}_0=s(D,\bold P^3)_0+(a+b)H\cap s(D,\bold P^3)_1
\label{eq:fulton1}
\end{equation}
where $c(N)$ denotes the
total Chern class of $N$, and $s(D,\bold P^3)$
is the Segre class of $D$ in $\bold P^3$ (\cite{fu}, $\S$9 ).
Summing up, we have
\begin{equation}
\tilde D\cdot E_C\le s(D,\bold P^3)_0+(a+b)H\cap s(D,\bold P^3)_1
\label{eq:fulton2}
\end{equation}
and equality holds if and only if
$D$ does not meet any component of $G-C$ different from $C$.

\begin{lem}
We have $s(D,\bold P^3)_1=[D]$ and
$s(D,\bold P^3)_0\le -2deg(D)$; if $C$ only
has ordinary singularities then equality holds
if and only if $D$ is a line.
\end{lem}
{\it Proof.}
If either
(a) or (b) in the statement of Lemma X holds, then
$D$ is a local complete intersection, and therefore it has
a normal bundle $N$ in $\bold P^3$.
Hence $s(D,\bold P^3)=c(N)^{-1}\cap [D]$, and
the statement is then reduced to the inequality
$deg(N)\ge 2deg(D)$.   Let $t:D_n@>>>D$ be the normalization
of $D$.
We then have a generically surjective morphism
$t^{*}T_{\bold P^3}@>>>t^{*}N$. On the other hand,
$T_{\bold P^3}(-1)$ is globally generated, and therefore we
must have $deg(N(-1))\ge 0$,
i.e. $deg(N) \ge 2d$.
If furthermore $D$ only has ordinary singularities,
we have an exact sequence $0@>>>T_D@>>>t^{*}T_{\bold P^3}
@>>>N$ and this shows that
equality holds if and only if
$g=0$ and $d=1$.
$\sharp$
\bigskip

We then have
$\tilde D\cdot E_C\le (a+b-2)deg(D)$, and if $D$ has only ordinary
singularities then
equality holds if and only if $D$ is a line not meeting any
component of $G-D$ different from $C$.
The Proposition follows.
$\sharp$
\bigskip

We know define two auxiliary invariants related to the Seshadri
constant that will be useful shortly.
\begin{defn} Let $C\subset \bold P^3$
be a smooth curve of degree $d$ and let $\epsilon (C)$ be its
Seshadri constant.
Let $N$ be the normal bundle of $C$ in $\bold P^3$.
For $0\le \eta \le \epsilon (C)$ a rational number,
define
$$\delta _{\eta}(C)=:\eta \cdot deg(N)-d$$
and
$$\lambda _{\eta}(C)=:\eta ^2d^2-\delta _{\eta}(C).$$
\label{defn:dandl}
\end{defn}
It is easy to check that
\begin{equation}
\delta _{\eta}(C)=:E^2\cdot H_{\eta}.
\label{eq:dint}
\end{equation}
More explicitly,
suppose that $0<\eta <\epsilon (C)$
and let $m$ and $n$ be large positive integers
such that $\eta =\frac nm$ and
$mH-nE$ is very ample.
Then for a general
$S\in |mH-nE|$ the intersection $C^{\prime}=E\cap S$ is
an irreducible smooth curve, and the induced morphism
$C^{\prime}@>>>C$ has degree $n$.
Then
\begin{equation}
\delta _{\eta}(C)=\frac {C^{\prime}\cdot _SC^{\prime}}
{H\cdot _SH}.
\label{eq:dint1}
\end{equation}
Similarly,
\begin{equation}
\lambda _{\eta}(C)=\frac {(H\cdot C^{\prime})^2}{(H\cdot _SH)^2}-
\frac {C^{\prime}\cdot _SC^{\prime}}
{H\cdot _SH}.
\label{eq:lint}
\end{equation}
\begin{rem} If we let $x=\eta d$, we have
$\lambda _{\eta}(C)=f(x)$, where
$$f(x)=x^2-\Big (4+\frac {2g-2}d\Big )x+d.$$
For $C$ subcanonical, $f$ is the polynomial introduced
by Halphen in his celebrated {\it speciality theorem}
(\cite{gp}), given by
$$g(X)=x^2-(4+e)x+d$$
where $e=max\{k|H^1(C,\cal O_C(k))\neq 0\}.$
Observe that $e\le (2g-2)/d$ always.
\end{rem}
\begin{cor}
Suppose that there exists an irreducible surface of degree $m$
through $C$, having multiplicity $n$ along $C$. If $\eta =\frac nm$,
then $\lambda _{\eta}(C)\ge 0$.
In particular, $\lambda _{\eta}(C)\ge 0$
for all $0\le \eta \le \epsilon (C)$.
Equality holds if and only if
$\cal O_S(C^{\prime})$ is numerically equivalent
to a multiple of $\cal O_S(H)$.
In particular, $\lambda _{\epsilon (C)}(C)\ge 0$
and equality holds if
$C$ is a complete intersection.
If $C$ is subcanonical and $\eta d$ is an integer, then
$\lambda _{\eta}(C)=0$
forces $C$ to be a complete intersection.
\label{cor:lpos}
\end{cor}
{\it Proof.} A straightforward application of the Hodge
index theorem.
The last part follows from the corresponding
statement of the speciality theorem (see \cite{gp}).
$\sharp$
\bigskip

\begin{cor} We have:
$$g\le \frac 12 d^2\epsilon (C)+d\Big (\frac 1{2\epsilon (C)}-2
\Big )+1.$$
\end{cor}

The right-hand side of the above inequality is
a decreasing function of $\epsilon$ in the interval
$(1/d, 1/{\sqrt d})$. In other words, higher Seshadri constants
impose tighter conditiond on the genus.
For a Castelnuovo extremal curve of even degree
we have $\epsilon =\frac 2d$ and the right hand side, as
a function of $d$, is asimptotic to $\frac {d^2}4$.

\begin{cor} Let $D$ be a divisor on $X_C$, and
set $s=D\cdot H_{\eta}\cdot H$.
Then for $0\le \eta \le \epsilon (C)$ we have
$$D^2\cdot H_{\eta}-D\cdot H_{\eta}\cdot E\le s^2-s\eta d.$$
\label{cor:sl}
\end{cor}

{\it Proof.} Write
$$D=xH+yE.$$
Then
$$D^2\cdot H_{\eta}=x^2+y^2\delta _{\eta}(C)+2xyd$$
and
$$D\cdot H_{\eta}\cdot E=x\eta d+y\delta _{\eta}(C).$$
{}From this we obtain
$$D^2\cdot H_{\eta}-D\cdot H_{\eta}\cdot E
=s^2-s\eta d-\lambda _{\eta}(C)(y^2-y).$$
Since $y$ is an integer, the statement then follows from
Corollary \ref{cor:lpos}.
$\sharp$
\bigskip

\begin{rem} From the inequality (see Remark \ref{rem:wahl})
$$s(N)\ge \frac 12 deg(N)$$
and the definition of $\delta _{\eta}(C)$, it
is easy to see that
$$d\ge \delta _{\eta}(C).$$
\label{rem:dd}
\end{rem}

\section{\bf {Gonality of space curves
and free pencils on projective
varieties}}\label{section:gon}

We have seen that if $C\subset S$ is a smooth curve with
$C^2>0$, then one can give lower bounds on the
gonality of $C$.
We deal here with the next natural question: if $C\subset \bold P^3$,
what can be said about $gon(C)$
in terms of the invariants of this embedding, and exactly which invariants
should one expect to play a direct role?
A hint to this is given by Lazarsfeld's result, to the
effect that if $C$ is nondegenerate complete intersection
of type $(a,b)$ with $a\ge b$  then
$gon(C)\ge a(b-1)$.

For $C\subset \bold P^r$ a smooth curve,
we let
$$\delta _{\eta}(C)=E^2\cdot H_{\eta}^2.$$
We then have
$\delta _{\eta}(C)=\eta ^{r-3}(\eta deg(N)-deg(C))$.

\begin{thm} Let $C\subset \bold P^r$ be a
smooth curve of degree $d$, $r\ge 3$.
Let $\epsilon (C)$ be the Seshadri constant of $C$, and
set $\alpha =min \Big \{1,\sqrt {\epsilon (C)
^{r-3}d}\Big (1-\epsilon (C)\sqrt {\epsilon (C)^{r-3} d}\Big
)\Big \}$.
Then
$$gon(C)\ge min\Big \{ \frac {\delta _{\epsilon (C)}(C)}{4\epsilon (C)
^{r-2}},
\alpha \Big (deg(C)-\frac {\alpha}{\epsilon (C)^{r-2}}\Big )\Big \}.$$
\end{thm}

Although we state the result for curves in $\bold
P^r$ for the sake of
simplicity, it is easy to see that
the same considerations apply when $\bold P^r$ is replaced by
a general smooth projective manifold $X$ with
$Pic(X)\simeq \bold Z$. Later in this section we shall indicate
how these results generalize to higher dimensional varieties
in $\bold P^r$.

{\it Proof.} To avoid
cumbersome notation, we shall assume that $r=3$. The proof
applies to higher value of $r$, with no significant change.
We want then to show that
\begin{equation}
gon(C)\ge min \Big \{\frac {\delta _{\epsilon (C)}(C)}
{4\epsilon (C)}, \alpha \Big (d-\frac {\alpha}{\epsilon (C)}\Big )
\Big \},
\label{eq:spcv}
\end{equation}
where $\alpha =min \{1,\sqrt d(1-\epsilon (C)\sqrt d)\}$.

Suppose, to the
contrary, that the statement is false: if $k=gon(C)
$, then $k$ is strictly smaller than both terms within the braces in the
last inequality. For $\eta <\epsilon (C)$ sufficiently close
to $\epsilon (C)$ the same inequality holds. More precisely,
if let $\alpha _{\eta}=min \{1,\sqrt d(1-\eta \sqrt d)\}$, we
have:
\begin{equation}
k<\dfrac {\delta _{\eta}(C)}{4\eta}
\label{eq:in1}
\end{equation}
and
\begin{equation}
k< \alpha _{\eta}\Big (d-\dfrac {\alpha _{\eta}}{\eta}\Big ).
\label{eq:in2}
\end{equation}
Pick a minimal pencil $A\in Pic^k(C)$, and set $V=:H^0(C,
A)$. Then $V$ is a two-dimensional vector space.
On $C$ we have an exact sequence of locally free sheaves $0@>>>
-A@>>>V\otimes \cal O_C@>>>A@>>>0$.
Consider the blow up diagram:
\begin{equation}
\CD
E @>>> X_C=Bl_C(X) \\
@V\pi VV     @VVfV \\
C@>>>    X
\endCD
\label{eq:cd}
\end{equation}
(here $E$ clearly denotes the exceptional divisor).
Define
\begin{equation}
\cal F=: Ker(\psi :V\otimes \cal O_{X_C}@>>>\pi ^{*}A).
\label{eq:F}
\end{equation}
$\pi ^{*}A$ is a line bundle on $E$, and $\psi$ is surjective.
Since $E$ is a Cartier divisor in $X_C$, $\cal F$ is a rank two
vector bundle on $X_C$.
As usual we set $H_{\eta}=H-\eta E$, where $\eta$ is a rational
number.
\begin{claim} Let $\eta$ be a rational number in the interval
$(0,\epsilon (C))$.
If $k<\dfrac {\delta _{\eta}(C)}{4\eta}$, then
$\cal F$ is Bogomolov-unstable with respect to $H_{\eta}$.
\end{claim}
{\it Proof.}
By Lemma \ref{lem:eltr},
the Chern classes of $\cal F$
are $c_1(\cal F)=-E$ and
$c_2(\cal F)=\pi ^{*}[A]$, where $[A]$ denotes the divisor class
in $A^1(C)$ of any element in $|V|$, and we implicitly map
$A^1(E)$ to $A^2(X_C)$.
Then the discriminant of $\cal F$ (definition \ref{defn:discr})
is given by
$$ \Delta (\cal F)=E^2-4[A].$$
Therefore by the assumption we have
\begin{equation}
\Delta (\cal F)\cdot H_{\eta}=
\delta _{\eta}(C)-4\eta k>0,
\label{eq:basin2}
\end{equation}
which implies that $\cal F$ is Bogomolov-unstable
with respect to $H_{\eta}$.
$\sharp$
\bigskip

Therefore, by Theorem \ref{thm:main}, there exists a
unique saturated invertible subsheaf
$\cal L\subset \cal F$ satisfying the following
properties:

(i) $\cal L$ is the maximal destabilizing subsheaf of $\cal F$ with
respect to any pair $(H_{\eta},R)$, with $R$ an arbitrary ample
divisor on $X_C$. In particular, for any such pair we have:
$(2c_1(\cal L)-c_1(\cal F))\cdot H_{\eta}\cdot R>0$.
Incidentally, this implies that $\cal L$ is the same for all
the values of $0< \eta <\epsilon (C)$ which make the hypothesis
of the claim true.

(ii) $(2c_1(\cal L)-c_1(\cal F))^2\cdot H_{\eta}\ge
\Delta (\cal F)\cdot H_{\eta}$.

\bigskip
Given the inclusions $\cal L\subset \cal F\subset \cal O_{X_C}^2$,
we have
\begin{equation}
\cal L=\cal O_{X_C}(-D)
\label{eq:eff}
\end{equation}
for some effective divisor $D$ on
$X_C$.
We can write
$$D=xH+yE,$$
with $x$ and $y$ integers and $x\ge 0$.
Set $$s=:D\cdot H_{\eta}\cdot H=x+y\eta d.$$
Since $\cal F$ has no sections, $D\neq 0$. The same applies for
the restriction to any ample surface. Hence $s\ge 0$
for $0<\eta <\epsilon (C)$.

\begin{lem} Assume that
$s\ge \alpha$. Then $k\ge \alpha(d-\dfrac {\alpha} {\eta})$.
\end{lem}
{\it Proof.} Given
(\ref{eq:eff}), from (ii) and (\ref{eq:basin2}) we get
\begin{equation}
(E-2D)^2\cdot H_{\eta}\ge \delta _{\eta}(C)-4\eta k.
\label{eq:basin3}
\end{equation}
Since $E^2\cdot H_{\eta}=\delta _{\eta}(C)$, this can be rewritten
$$D^2\cdot H_{\eta}-D\cdot H_{\eta}\cdot E\ge
-\eta k.$$
By Corollary \ref{cor:sl}, we then have
\begin{equation}
s^2-s\eta d\ge -\eta k.
\label{eq:basin4}
\end{equation}
On the other hand,
we have the destabilizing condition (i)
\begin{equation}
(E-2D)\cdot H_{\eta}\cdot H\ge 0.
\label{eq:dest}
\end{equation}
Now
$$E\cdot H_{\eta}\cdot H=\eta d$$
and therefore (\ref{eq:dest}) can be written
\begin{equation}
\eta d\ge 2s.
\label{eq:dest1}
\end{equation}
Therefore we can apply Lemma \ref{lem:trivial} with $a=\eta d$ and
$b=\eta k$ to obtain
$$\eta k\ge \eta d \alpha-\alpha ^2.$$
This proves the Lemma.
$\sharp$
\bigskip

The proof of the theorem is then reduced to the following Lemma.
\begin{lem} $s\ge \alpha$.
\end{lem}
{\it Proof.} We shall argue that
$s\ge \alpha _{\eta}$ for all rational $\eta <\epsilon (C)$ such
that the inequalities (\ref{eq:in1}) and
(\ref{eq:in2}) hold.
For all such $\eta$ we are then in the
situation of Claim 4.1.

\begin{claim} $\cal L$ is saturated in $V\otimes \cal O_X$.
\label{claim:sat}
\end{claim}
{\it Proof.} By construction, $\cal L
=\cal O_X(-D)$ is saturated in $\cal F$.
Therefore, if the Claim is false then the inclusion
$\cal L\subset V\otimes \cal O_X$ drops rank along $E$.
Hence there exists an inclusion
$\cal O_X(E-D) \subset \cal O_X^2$.
This implies that $D-E$ is effective, and in particular
$(D-E)\cdot H_{\eta}^2\ge 0$.
Together with
the instability condition $(E-2D)\cdot H_{\eta}^2>0$,
this would imply
$D\cdot H_{\eta}^2<0$,
against the fact that $D$ is effective.
$\sharp$
\bigskip

By Claim \ref{claim:sat}, there is an exact sequence
$$0@>>>\cal O_X(-D)@>>>V\otimes \cal O_X@>>>\cal O_X(D)\otimes
\cal J_Y@>>>0,$$
where $Y$ is a closed subscheme of $X$ of codimension two or empty.
Computing $c_2(\cal O_X^2)=0$ from this sequence, we obtain
$D^2=[Y],$
and therefore
$D^2\cdot H\ge 0$.
On the other hand,
$D^2\cdot H=x^2-y^2d$, and so
$$x\ge |y|\sqrt d.$$
Now,
$$s=x+y\eta d\ge x-|y|\eta d\ge |y|\sqrt d(1-\eta \sqrt d).$$
By construction, $H^0(X,\cal F)=0$, and therefore
$D\neq 0$. Hence, if $y=0$ then $s=x\ge 1$.
If $y\neq 0$, then the above inequality shows
that
$s\ge \sqrt d(1-\eta \sqrt d)$.
$\sharp$
\bigskip

This completes the proof of the Theorem.
$\sharp$
\bigskip

\begin{exmp} Let $C\subset \bold P^3$ be a smooth complete
intersection curve of type $(a,b)$, with $a\ge b+3$, $b\ge 2$.
Then $gon(C)\ge a(b-1)$.
\label{exmp:ci}
\end{exmp}
\begin{rem}
This shows that
the result is generally optimal.
However,
the theorem is void for a complete intersection of type $(a,a)$.
But for complete intersections
one knows more than just the Seshadri constant: not
only $\epsilon (C)=\frac 1a$, but
in fact the linear series $|aH-E|$ is base point free, and the
general element is smooth.
An ad hoc argument proves that $gon(C)\ge a(b-1)$ (\cite{la:unp}).
\label{rem:ci}
\end{rem}

\begin{exmp} Let $C$ be a nondegenerate smooth complete
curve in $\bold P^3$ that is linked to a line in a
complete intersection of type $(a,b)$. Then for
$a\gg b\gg 0$ we obtain
$gon (C)\ge deg(C)-(a+b-2)$.
This is clearly optimal,
because a base point free linear series of that degree
is obtained by considering the pencil of planes through the
line.
The same considerations as in Remark
\ref{rem:ci} apply.
\end{exmp}

\begin{rem} An analysis of "small" linear series on special
classes of space curves is carried out by Ciliberto and
Lazarsfeld in \cite{cl}. It would be interesting to know whether the
present method can be adapted to give a generalization of
their results.
\end{rem}

{}From the Theorem, we immediately get
\begin{cor} Let $X\subset \bold P^r$ be
a smooth projective variety.
Let
$d$ be the degree of $X$ and $\epsilon (X)$ be its Seshadri constant.
Suppose that $A$ is a line bundle on $X$ with a pencil
of sections $V\subset H^0(X,A)$ whose base locus has codimension
at least two. Let $F$ be any divisor in the linear series
$|A|$. Then
$$deg(F)\ge
min \Big \{
\frac 1{4\epsilon (C)^{r-2}} \Big [\epsilon (X) \big (c_1(N)\cdot _XH^{n-1}+
(n-1)d \big )-d\Big ],
\alpha \Big (d-\frac {\alpha}{\epsilon (X)}\Big )\Big \},$$
where
$\alpha =\Big \{1,\sqrt {\eta ^{r-3}d}
\Big (1-\epsilon (X)\sqrt {d\epsilon (C)^{r-3}} \Big )\Big \}$.
\end{cor}
{\it Proof.} Let $C\subset X$ be a curve
of the form $X\cap \Lambda$, where
$\Lambda\subset \bold P^3$ is a linear subspace of dimension
$c+1$, with $c$ the codimension of $X$.
Then $V$ restricts to a base point free pencil on $C$, and
the result follows by applying the theorem.
$\sharp$
\bigskip

Given the general nature of the above arguments,
one clearly expects that they should be applicable to a
wider range of situations. In fact,
we give now the generalization of theorem
3.1 to arbitrary smooth projective varieties in $\bold P^r$.
The proof is exactly the same as the one for theorem 3.1,
the only change consisting in a more involved notation.

\begin{thm} Let $Y\subset \bold P^r$ be a projective manifold
of degree $d$ and codimension $c$. Let $\epsilon (Y)$ be its Seshadri
constant, and suppose $0\le \eta \le \epsilon (Y)$.
If $A$ is base point free pencil on $Y$, then
$$\pi ^{*}[A]\cdot H_{\eta}^{r-2}\ge
min \Big \{\frac {\delta _{\eta}(Y)}4, \frac 1{H^2H_{\eta}^{r-2}}
\alpha \Big (H\cdot E\cdot H_{\eta}^{r-2}- \alpha \Big ) \Big \},$$
where
$\alpha =min \Big \{H^2\cdot H_{\eta}
^{r-2}, \sqrt {\eta ^{c-2} deg(Y)}
\Big (H^2\cdot H_{\eta}^{r-2}-\dfrac {H\cdot E\cdot H_{\eta}^{r-2}}
{\sqrt {\eta ^{c-2} deg(Y)}}\Big ) \Big \}$.
\end{thm}

\section{\bf {Stability of restricted bundles}}\label{section:stability}

We deal here with the following problem:
\begin{prob}
Let $\cal E$ be a rank two vector bundle on $\bold P^3$,
and let $C\subset \bold P^3$ be a smooth curve.
If $\cal E$ is stable, what conditions on $C$ ensure
that $\cal E|_C$ is also stable?
\end{prob}

\begin{rem}
This question has been considered by Bogomolov
(\cite{bo:78} and \cite{bo:svb})
in the case of vector bundles on surfaces.
In particular, he shows that if $S$ is a smooth projective
surface with $Pic(S)\simeq \bold Z$, $\cal E$ is a stable rank
two vector bundle on $S$ with $c_1(\cal E)=0$ and $C\subset S$
is a smooth curve with $C^2>4c_2(\cal E)^2$, then
$\cal E|_C$ is stable.
\label{rem:surfacecase}
\end{rem}

After a suitable twisting, we may also assume that
$\cal E$ is {\it normalized}, i.e.
$c_1(\cal E)=0$ or $-1$. We shall suppose here
that $c_1(\cal E)=0$, the other case being similar.

As usual we adopt the following notation:
$$f:X_C@>>>\bold P^3$$
is the blow up of $\bold P^3$ along $C$,
$$E=f^{-1}C$$
is the exceptional divisor, and
$$\pi :E@>>>C$$
is the induced projection.
Recall that for $\eta \in \bold Q$ we set
$$H_{\eta}=:H-\eta E,$$
where we write $H$ for $f^{*}H$.
If $0<\eta <\epsilon (C)$, $H_{\eta}$ is a polarization
on $X_C$.
\begin{defn} We define the {\it stability constant}
of $\cal E$ w.r.t. $C$ as
\begin{center}
$\gamma (C,\cal E)=
sup \{\eta \in [0,\epsilon (C)]| f^{*}\cal E$ is $(H_{\eta},
 H)$-stable$\}$.
\end{center}
\label{defn:gamma}
\end{defn}
\begin{rem}
Recall that $f^{*}\cal E$ is $(H,H_{\eta})$-stable if for all
line bundles $\cal L\subset f^{*}\cal E$ we have
$\cal L\cdot H\cdot H_{\eta}<0$.
\end{rem}

\begin{lem} Suppose $0\le \eta <\epsilon (C)$.
Then $f^{*}\cal E$ is $(H,H_{\eta})$-semistable if and only if
$\eta \le \gamma (C,\cal E)$. If
$\eta <\gamma (C,\cal E)$, $f^{*}\cal E$ is
$(H,H_{\eta})$-stable.
\end{lem}
{\it Proof.} The collection of the numerical classes
of nef divisors $D$ with respect to which $f^{*}\cal E$
is $(H,D)$-semistable (or stable) is convex, hence it
contains
the segment $[H,H_{\gamma (C,\cal E)}]$.
Since $f^{*}\cal E$ is $(H,H)$-stable the second
statement follows.
$\sharp$
\bigskip

\begin{lem} Suppose that $V\subset \bold P^3$ is a smooth
surface of degree $a$ containing $C$, and that
$\cal E|_V$ is $\cal O_V(H)$-stable.
Then
$$\gamma (C,\cal E)\ge min\big \{\frac 1a,\epsilon (C)\big \}.$$
\label{lem:bound}
\end{lem}
{\it Proof.} Let $\tilde V$ be the proper transform
of $V$ in $X_C$.
We have $\tilde V\simeq V$
and
$$\tilde V\in |aH_{\frac 1a}|.$$
The hypothesis implies that for every line-bundle
$$\cal L\subset  f^{*}\cal E$$
we have
$$\cal L\cdot H_{\frac 1a}\cdot H<0.$$
Hence the same holds for every
$\eta$ with $0\le  \eta \le \frac 1a$.
$\sharp$
\bigskip

\begin{rem} Note that the same argument actually proves the
following stronger statement: let $V\supset C$ be a reduced
irreducible surface
through $C$ having degree $m$ and multiplicity $n$ along $C$,
and such that $f^{*}\cal E|_{\tilde V}$ is $\cal O_{\tilde V}(H)$-stable.
Then $\gamma (C,\cal E)\ge min\{\frac nm,\epsilon (C)\}$.
\end{rem}

\begin{lem} Fix $c_2\ge 0$ an integer.
Then then there exists an integer $k$ with the
following property. If $\cal E$ is a
stable rank two vector bundle
on $\bold P^3$ with $c_1(\cal E)=0$ and $c_2(\cal E)=c_2$,
and if $V\subset \bold P^3$ is a smooth surface of
degree $a\ge k$, then $
\cal E|_V$ is $\cal O_V(H)$-stable.
\label{lem:surfres}
\end{lem}
{\it Proof.} We start by finding $s$ such that for a general
surface $S$ of degree $s$
we have $Pic(S)\simeq \bold ZH$ ($s\ge 4$ will do) and furthermore
the restriction $\cal E|_S$ is $\cal O_S
(H)$-stable.
Bogomolov's theorem
(remark \ref{rem:surfacecase}) then says that for any
curve $C\subset S$ such that $C^2>4c_2(\cal E)^2s^2$ the
restriction $\cal E|_C$ is also stable.
Let now $a>0$ be such that $a^2>4c_2(\cal E)^2s$.
Suppose that $V$ is a smooth surface of degree $a$ and that
$\cal E|_V$ is not stable. Then the same is true for
$C=V\cap S$. For a general choice of $S$, $C$ is a smooth
curve, and since $C\cdot _SC=a^2s>4c_2(\cal E)^2s^2$,
we have a contradiction.
$\sharp$
\bigskip

We can in fact restate the previous lemma as follows:

{\it Let $s$ be the smallest positive integer such that
for a general surface of degree $s$ we have $Pic(S)\simeq \bold
Z$ and $\cal E|_S$ stable.
If $a>2c_2(\cal E)\sqrt s$ and $V\subset \bold P^3$ is any smooth surface
of degree $a$, then $\cal E|_V$ is $\cal O_V(H)$-stable.
}

\begin{cor} Let $\cal E$ be a rank two stable bundle on $\bold P^3$
with $c_1(\cal E)
=0$ but $c_2(\cal E)\neq 1$ (i.e., $\cal E$ is not a null correlation
bundle (\cite{oss}). If $a>2c_2(\cal E)$ and $V\subset \bold P^3$
is a smooth hypersurface of degree $a$, then $\cal E|_V$ is
$\cal O_V(H)$-stable.
\label{cor:surfres}
\end{cor}
{\it Proof.} In fact, a theorem of Barth says that in this case
we can take $s=1$ (\cite{ba}).
$\sharp$
\bigskip

\begin{rem} In light of Barth's restriction theorem,
by induction these statements
generalize to $\bold P^r$ for any $r\ge 2$ (for $r=2$ this is just
Bogomolov's theorem, and the hypothesis $c_2\neq 1$ is not needed).
\end{rem}

\begin{rem}
In the proof of Corollary 5.1, we use stability on the
hyperplane section to deduce stability on the whole surface.
What makes this work is Bogomolov's theorem
(cfr remark \ref{rem:surfacecase}), which gives us
a control of the behaviour of stability under restriction to
plane curves. On the other hand, if we are given an arbitrary
smooth surface $V$, it may well be that $\cal E|_V$ is $H$-stable
while $\cal E|_C$ is not, where $C$ is an hyperplane section of
$V$. In that case, however, $\cal E|_{V\cap W}$ will be stable,
if $W$ is a smooth surface of
very large degree such that $V\cap W$
is a smooth curve. To improve the above result, therefore,
one would need to control the behaviour of stability under restriction
to curves not necessarily lying in a plane.
After proving the restriction theorem \ref{thm:curveres}
we shall strengthen
the above corollary (cfr Corollary \ref{cor:surres1}).
\end{rem}
\bigskip

\begin{defn}
If $X$ is a smooth variety and
$c_i\in A^i(X)$ for $i=1$ and $2$,
let $\cal M_X(c_1,c_2)$ denote the
moduli space of stable rank two vector bundles with Chern classes
$c_1$ and $c_2$.
\end{defn}

\begin{cor} Fix an integer $c_2\ge 0$. Then for
any sufficiently large positive integer $a$ the following
holds: if $V\subset \bold P^r$ is a smooth hypersurface of degree
$a$, then $\cal M_{\bold P^r}(0,c_2)$ embeds as an open subset of
$\cal M_V(0,c_2a)$.
\label{cor:surfmoduli}
\end{cor}
{\it Proof.} $\cal M_{\bold P^r}(0,c_2)$ forms a bounded
family of vector bundles, and therefore so does the collection
of the vector bundles $End (\cal E,\cal F)$, with
$\cal E$, $\cal F\in \cal M_{\bold P^r}(0,c_2)$.
Therefore, if $k\gg 0$, we have
$$H^i(\bold P^r,End (\cal E,\cal F)(-a))=0$$
for all
$i\le 2$,
$a\ge k$ and for all $\cal E$, $\cal F\in \cal M_{\bold
P^r}(0,c_2)$. Furthermore, by the above
lemma we can assume that $\cal E|_V$ is
$\cal O_V(H)$-stable for all
$\cal E\in \cal M_{\bold P^r}(0,c_2)$.
{}From the long exact sequence
in cohomology associated to the exact sequence of sheaves
$$0@>>>End(\cal E,\cal F)(-a)@>>>End(\cal E,\cal F)@>>>
End(\cal E|_V,\cal F|_V)@>>>0$$
we then obtain
isomorphisms
$$H^0(\bold P^r,End(\cal E,\cal F))\simeq H^0(V,
End (\cal E|_V,\cal F|_V))$$
and
$$H^1(\bold P^r,End (\cal E,\cal F))\simeq
H^1(V,End (\cal E|_V,\cal F|_V)).$$
Since there can't be any
nontrivial homomorphism between two nonisomorphic
stable bundles
of the same slope, the first one says that
$$\cal E@>>>\cal E|_V$$
is a one-to-one morphism of $\cal M_{\bold P^r}(0,c_2)$ into
$\cal M_V(0,c_2a)$ and the second says that
the derivative of this morphism is an isomorphism
(\cite{ma}).
$\sharp$
\bigskip

\begin{cor} $\gamma (C,\cal E)>0$.
\label{cor:gammapos}
\end{cor}
{\it Proof.} By Lemma 5.3, for $r\gg 0$
the restiction of $\cal E$ to any smooth surface of degree $r$
is stable with respect to the hyperplane bundle. So we
just need to consider a smooth surface through $C$ of very large
degree and apply Lemma 5.2.
$\sharp$
\bigskip

\begin{exmp} Let
$$C=V_a\cap V_b\subset \bold P^3$$
be
a smooth complete intersection of type $(a,b)$, with
$a\ge b$.
Suppose that $V_a$ is smooth, and that
$\cal E|_{V_a}$ is $\cal O_{V_a}(H)$-stable.
Then
$$\gamma (C,\cal E)=\frac 1a=\epsilon (C).$$
In general, $0<\eta <\gamma (C,\cal E)$ if
and only if for $m$ and $n$ sufficiently large
integers such that $\eta =\frac nm$, and $S\in |mH-nE|$
a smooth surface, we have
that $f^{*}\cal E|_S$ is $\cal O_S(H)$-stable. In other words,
we have a degree $m$ hypersurface with an ordinary singularity
of multiplicity $n$ along $C$, such that the pull-back
of $\cal E$ to the desingularization of $S$ is $H$-stable.
\end{exmp}
\bigskip
Our result is then the following:
\begin{thm} Let $\cal E$ be a stable rank two vector
bundle on $\bold P^3$ with $c_1(\cal E)=0$.
Let $C\subset \bold P^3$ be a smooth curve
of degree $d$ and Seshadri constant $\epsilon (C)$,
and let $\gamma =\gamma (C,\cal E)$ be the stability constant
of $\cal E$ w.r.t $C$.
Let
$\alpha =min\Big \{1,\sqrt d\Big (\sqrt {\frac 34}-\gamma \sqrt d
\Big )\Big \}.$
Suppose that
$\cal E|_C$ is not stable.
Then
$$c_2(\cal E)\ge min\Big \{\frac {\delta _{\gamma}}4,
\alpha \gamma  \Big (d-\frac {\alpha}{\gamma}\Big )\Big \}.$$
\label{thm:curveres}
\end{thm}

{\it Proof.} Suppose to the contrary that $c_2(\cal E)$ is strictly
smaller that both quantities on the right hand side. We can
find a rational number $\eta$ with $0<\eta <\gamma $
such that
\begin{equation}
c_2(\cal E)<\dfrac {\delta _{\eta}(C)}4
\label{eq:c2ineq}
\end{equation}
and
\begin{equation}
c_2(\cal E)<\alpha \eta \Big (d-\dfrac {\alpha}{\eta}\Big).
\label{eq:c2ineq1}
\end{equation}

By assumption there exists a line bundle $L$ on $C$ of
degree $l\ge 0$ sitting in an exact sequence
$$0@>>>L@>>>\cal E|_C@>>>L^{-1}@>>>0.$$
Define a sheaf $\cal F$ on $X_C$ by the exactness of the
sequence
\begin{equation}
0@>>>\cal F@>>>f^{*}\cal E@>>>\pi ^{*}L^{-1}@>>>0.
\label{eq:Fi}
\end{equation}
Then $\cal F$ is a rank two vector bundle on $X_C$
having Chern classes
$c_1(\cal F)=-[E]$ and $c_2(\cal F)=f^{*}c_2(\cal E)-\pi
^{*}[L]$ (cfr Lemma \ref{lem:eltr}).
A straightforward computation then gives
\begin{equation}
\Delta (\cal F)\cdot H_{\eta}=
\delta _{\eta}(C)-4c_2(\cal E)+4\eta l\ge \delta _{\eta}(C)-4c_2(\cal E)
\label{eq:DF}
\end{equation}
and this is positive by (\ref{eq:c2ineq}).
Therefore $\cal F$ is Bogomolov-unstable with respect
to $H_{\eta}$ (Theorem \ref{thm:main}). Let
$$\cal L\subset \cal F$$
be the maximal destabilizing line bundle w.r.t. $H_{\eta}$. We
shall write
$$\cal L=\cal O_{X_C}(-D),$$
with
$$D=xH-yE.$$
\begin{claim}  $x>0$
\label{claim:xpos}
\end{claim}
{\it Proof.} Pushing forward the inclusion
$\cal L\subset \cal F$ we obtain an inclusion
$\cal O_{\bold P^3}(-x)\subset \cal E$. Therefore the
statement follows from the assumption of stability on
$\cal E$ and the hypothesis $c_1(\cal E)=0$.
$\sharp$
\bigskip

The destabilizing condition says
$(2c_1(\cal L)-c_1(\cal F))\cdot H_{\eta}\cdot R\ge 0$
for all nef divisors on $X_C$, with strict inequality holding
when $R$ is ample. In particular, with $R=H$ we have
\begin{equation}
(E-2D)\cdot H_{\eta}\cdot H\ge 0.
\label{eq:resinst}
\end{equation}
Let us set $s=D\cdot H_{\eta}\cdot H$. Then (\ref{eq:resinst})
reads
\begin{equation}
\eta d\ge 2s.
\label{eq:resinst1}
\end{equation}
On the other hand, since $\cal L$ is saturated in
$\cal F$, we also have
$(E-2D)^2\cdot H_{\eta}\ge \Delta (\cal F)\cdot H_{\eta}$,
and with some algebra this becomes
\begin{equation}
c_2(\cal E)\ge D\cdot E\cdot H_{\eta}-D^2\cdot H_{\eta}+\eta l\ge
D\cdot E\cdot H_{\eta}-D^2\cdot H_{\eta}.
\label{eq:resinst2}
\end{equation}
Invoking Corollary \ref{cor:sl}, we
then have
\begin{equation}
c_2(\cal E)\ge s\eta d-s^2.
\label{eq:resinst3}
\end{equation}

\begin{claim} $\cal L$ is saturated in $f^{*}\cal E$.
\label{claim:ressat}
\end{claim}
{\it Proof.} Suppose not. Then there would be an inclusion
$$\cal L(E)=\cal O_{X_C}(E-D)\subset f^{*}\cal E$$
and therefore the $(H,H_{\eta})$-stability of $f^{*}\cal E$
would force
$$(E-D)\cdot H_{\eta}\cdot H<0.$$
On the other hand by instability we have $E\cdot H_{\eta}\cdot H\ge
2D\cdot H_{\eta}\cdot H$ and from this we would get
$$E\cdot H_{\eta}\cdot H=\eta d<0,$$
absurd.
$\sharp$
\bigskip

Therefore there is an exact sequence
$$0@>>>\cal O_{X_C}(-D)@>>>f^{*}\cal E@>>>\cal O_{X_C}(D)\otimes \cal J_W
@>>>0$$
where $W$ is a local complete intersection subscheme of $X_C$ of codimension
two or empty.
Computing $c_2(f^{*}\cal E)$ from the above
sequence we then get
$f^{*}c_2(\cal E)=W-D^2$, i.e.
$$D^2\cdot H\ge -c_2(\cal E).$$
This can be rewritten
$$x^2\ge y^2d-c_2(\cal E).$$
Recall that we have (remark \ref{rem:dd})
$$d\ge \delta _{\eta}(C),$$
and therefore the assumption $c_2(\cal E)< \delta _{\eta}(C)/4$
implies
\begin{equation}
c_2(\cal E)< \frac d4.
\label{eq:c2d4}
\end{equation}

\begin{lem}
$$s\ge min\Big \{1,\sqrt d\Big (\sqrt {\frac 34}-\eta \sqrt d\Big )
\Big \}.$$
\end{lem}
{\it Proof.}
If $y\le 0$ then $s=x+|y|\eta d\ge 1$.
If $y>  0$,
we have
$s=x-y\eta d\ge
y\sqrt d
\Big (\sqrt {1-\frac {c_2(\cal E)}d} -\eta\sqrt d\Big )$
and therefore using (\ref{eq:c2d4}) we obtain
$$s\ge
\sqrt d\Big (\sqrt {\frac 34}-\eta \sqrt d \Big ).$$
$\sharp$
\bigskip

Hence we can apply lemma
\ref{lem:trivial} with
$a=\eta d$ and $b=c_2(\cal E)$ to deduce
$$c_2(\cal E)\ge \alpha \eta d-\alpha ^2,$$
which contradicts (\ref{eq:c2ineq1}).
This completes the proof of the Theorem.
$\sharp$
\bigskip

\begin{cor} Let $\cal E$ be a stable rank two vector bundle on
$\bold P^3$ with $c_1(\cal E)=0$ and $c_2(\cal E)=c_2$.
If $b\ge c_2+2$ and $V\subset \bold P^3$ is a smooth hypersurface
of degree $b$, then $\cal E|_V$ is $\cal O_V(H)$-stable.
\label{cor:surres1}
\end{cor}
{\it Proof.}
Let $a\gg b$; then we may assume that if $W\subset \bold P^3$ is
a surface of degree $a$ then $\cal E|_W$ is $H$-stable.
If $W$ is chosen generally, we may also assume that $C=W\cap V$
is a smooth curve. Then
by Lemma \ref{lem:bound} we have $\gamma (C,\cal E)=\epsilon (C)=
\frac 1a$. For $a$ large enough, furthermore, we also have
$\alpha =1$. Hence the theorem says that if $\cal E|_C$ is not
stable, then $c_2\ge b-1$. The hypothesis implies therefore that
$\cal E|_C$ is stable, and this forces $\cal E|_V$ to be stable
also.
$\sharp$
\bigskip

\begin{cor}
Let $\cal E$ be as above, and let
$C=V_a\cap V_b$ be a smooth complete intersection
curve of type $(a,b)$, and suppose that $V_a$ is smooth.
Assume that $a\ge \frac 43b+\frac {10}3$ and that $b\ge c_2+2$.
Then $\cal E|_C$ is stable.
\end{cor}
{\it Proof.} By Corollary 5.4, $\cal E|_{V_b}$ is $H$-stable.
Hence by Lemma \ref{lem:bound} $\gamma (C,\cal E)=\frac 1a$.
The first hypothesis implies that $\alpha =1$. Hence if $\cal E|_C$
is not stable the theorem yields $c_2\ge b-1$, a contradiction.
$\sharp$
\bigskip

\begin{cor} Fix a nonnegative integer $c_2$. Then
we can find positive integers $a$ and $b$ such that if $C\subset
\bold P^3$ is any smooth complete intersection curve
of type $(a,b)$ then $\cal M_{\bold P^3}(0,c_2)$ embeds
as a subvariety of $\cal M_C(0)$.
\end{cor}
{\it Proof.} The argument is similar
to the one in the proof
of Corollary 5.2. Here one uses
the Koszul resolution
$$0@>>>\cal O_{\bold P^3}(-a-b)@>>>\cal O_{\bold P^3}(-a)\oplus
\cal O_{\bold P^3}(-b)@>>>\cal J_C@>>>0$$
to show that $H^i(\bold P^3,End(\cal E,\cal F)\otimes
\cal J_C)=0$ for $i\le 1$.
$\sharp$
\bigskip

\begin{rem} Using the above corollary,
we obtain a compactification of
$\cal M_{\bold P^3}(c_1,c_2)$, by simply
taking its closure in the moduli space of semistable bundles on the
curve.
It would be interesting
to know whether these compactifications are intrinsic, i.e. they
are independent of the choice of the curve or, if not,
how they depend on the geometry of the embedding $C\subset \bold P^3$.
\end{rem}


\begin{thebibliography}{Dillo 99}

\bibitem[ACGH]{acgh} Albarello, Cornalba, Griffiths and
Harris, {\em The geometry of Complex curves, I }Springer
Verlag 1985


\bibitem[Ba 77]{ba} W. Barth {\em Some properties of stable
rank two vector bundles on $\bold P^n$} Math. Ann. 226,
125-150 (1977)


\bibitem[BEL 89]{bel} A. Bertram, L. Ein and R. Lazarsfeld
{\em Vanishing theorems, a theorem of Severi and the equations defining
projective varieties} preprint 1989


\bibitem[Bo 78]{bo:78} F. Bogomolov {\em Unstable vector bundles
and curves on surfaces} Proc. ICM Helsinki (1978)



\bibitem[Bo 79]{bo:st} F. Bogomolov {\em Holomorphic tensors
and vector bundles on projective varieties} Math. USSR Izvestija
vol. 13 n. 3 (1979)


\bibitem[Bo 92] {bo:svb} F. Bogomolov {\em On stability of vector bundles
on surfaces and curves}
preprint


\bibitem[Ca 93] {ca} G. Castelnuovo {\em Sui multipli di una serie
lineare di punti appartenente a una curva algebrica}
Rend. Circ. Mat. di Palermo 7 (1893) 89-100



\bibitem[CL 84]{cl} C. Ciliberto and R. Lazarsfeld
{\em On the Uniqueness of certain linear series
on some classes of space curves} in Conte and Strano, editors,
L.N.M. 1092, Springer Verlag 1984


\bibitem[De 90] {de} J. P. Demailly {\em Singular hermitian metrics
on positive line bundles} Republications de l'Insitut Fourier
Gr\'enoble 1990



\bibitem[DM 89] {dm} R. Donagi and D. Morrison
{\em Linear systems on K3 sections} J. Diff. Geo. 28 (1989)
(49-64)


\bibitem[EL 92]{el} L. Ein and R. Lazarsfeld {\em Seshadri constants
on surfaces} preprint 1992




\bibitem[Fl 84] {fl} H. Flenner {\em Restrictions of
semistable bundles on projective varieties}
Comm. Maht. Elv. 59 (1984) 635-650


\bibitem[FHS 80]{fhs} O. Forster, A. Hirshowitz, M. Schneider
{\em Type de scindage g\'eneralis\'e pour les fibr\'es stables}
In: Vector Bundles and Differential Equations, p. 65-81
(Nice 1979) Birkhauser (1980)



\bibitem[Fu 84] {fu} W. Fulton {\em Intersection Theory}
Springer Verlag 1984


\bibitem[Gi 79] {gi} D. Gieseker {\em On a theorem of Bogomolov
on Chern classes of stable bundles} Am. J. Math. 101 (1979)
77-85



\bibitem[Gr 87]{gr} M. Green {\em Koszul Cohomology and
Geometry} in Cornalba {\it et al}, editors, Proceedings of
the ICTP College on Riemann Surfaces (1987), World Scientific
Press 1989


\bibitem[GL 87]{gl} M. Green and R. Lazarsfeld {\em Special divisors
on curves on a K3 surface} Inv. Math. 89 (1987) 635-650



\bibitem[GLP 83]{glp} L. Gruson, C. Peskine and R. Lazarsfeld
{\em On a theorem of Castelnuovo, and the equations defining space curves}
Inv. Math. 72 491-506 (1983)


\bibitem[GP 77] {gp}  L. Gruson and C. Peskine
{\em Th\'eoreme de specialit\'e}
S\'eminaire E.N.S. (1977-1978) Expos\'e XIII Axt\'erisque
71-72

\bibitem[Ha 77]{ha:ag} R. Hartshorne {\em Algebraic Geometry}
Springer Verlag 1977



\bibitem[La 87] {la:svbt} R. Lazarsfeld
{\em A sampling of vector bundle techniques in the study of
linear series}
in Cornalba {\it et al}, editors, Proceedings
of the ICTP College on Riemann Surfaces (1987)
World Scientific Press 1989

\bibitem[La]{la:unp} R. Lazarsfeld, unpublished


\bibitem[Ma 78] {ma} Maruyama {\em moduli of stable sheaves II}
J. of Kyoto Univ. 18 (557-614) (1978)


\bibitem[MR 82] {mr:res} V. B. Mehta and A. Ramanathan
{\em Semistable sheaves on projective varieties and
their restriction to curves}
Math. Ann. 258, 213-224 (1982)

\bibitem[MR 84] {mr:rep} V. B. Mehta and A. Ramanathan
{\em Restriction of stable shaves and representations of
the fundamental group} Inv. Math. (1984) (163-172)



\bibitem[Mi 85] {mi:cc} Y. Miyaoka
{\em The Chern classes and Kodaira dimension of a minimal variety}
Adv. Studies in pure Math. 10, (1987)
Alg. Geo., Sendai, 1985 (448-476)


\bibitem[Mu 66] {mu:cs} D. Mumford {\em Lectures on curves on an algebraic
surface} Ann. of Math. studies 59 (1966)


\bibitem[Mu 70]
{mu} D. Mumford {\em Varieties defined by quadratic
equations} in {\it Questions in algebraic varieties}, Centro Internazionale
di matematica estivo, Cremonese Roma 1970 (29-100)

\bibitem[OSS 80]{oss} C. Okonek, M. Schneider and
H. Spindler {\em Vector bundles on complex projective spaces}
Progress in Math. 3 Birkhauser 1980


\bibitem[Pa 93] {pa} R. Paoletti {\em Free pencils on divisors}
preprint


\bibitem[Re 88]{re:vbls} I. Reider
{\em Vector bundles of rank two and linear systems
on algebraic surfaces}
Ann. of Math. 127 (1988) 309-316


\bibitem[Re 89]{re:app} I. Reider {\em Applications of Bogomolov's
theorem} in Ciliberto {\it et al}, editors, Academic Press 1989

\bibitem[Se 87] {se:ext} F. Serrano {\em Extension of morphisms
defined on a divisor} Math. Ann. 277 (1987)

\bibitem[Sev 32]{sev}
F. Severi
{\em Le role de la g\'eometrie alg\'ebrique dans
les math\'ematiques}
Proc. Int. Cong. Math. (Zurich) 1932 vol. I 209-220

\bibitem[So 76] {so:amp} A. J. Sommese
{\em On manifolds that cannot be ample divisors}
Math. Ann. 1976 (221) (55-72)





\bibitem[Ty 87] {ty} A. N. Tyurin {\em Cycles, curves and vector
bundles on an algebraic surface} Duke Math. J. 54 (1987) (1-26)

\bibitem[Wa 91]{w} J. Wahl {\it Cohomology of symmetric bundles on
normal surface singularities} preprint 1991


\end{thebibliography}
\end{document}